\documentclass{article}%
\usepackage{graphicx}
\usepackage{amsmath}
\usepackage{amsfonts}
\usepackage{amssymb}%
\newtheorem{theorem}{Theorem}

\newtheorem{corollary}[theorem]{Corollary}

\newtheorem{definition}[theorem]{Definition}

\newtheorem{lemma}[theorem]{Lemma}

\begin{document}

\title{Gosset Polytopes in Picard groups \\of del Pezzo Surfaces}
\author{Jae-Hyouk Lee}
\maketitle

\begin{abstract}
In this article, we research on the correspondences between the geometry of
del Pezzo surfaces $S_{r}$ and the geometry of Gosset polytopes $(r-4)_{21}$.
We construct Gosset polytopes $(r-4)_{21}\ $in $Pic\ S_{r}\otimes\mathbb{Q}$
whose vertices are lines, and we identify divisor classes in $Pic\ S_{r}%
\ $corresponding to $(a-1)$-simplexes ($a\leq r$), $(r-1)$-simplexes and
$(r-1)$-crosspolytopes of the polytope $(r-4)_{21}$. Then we explain these
classes correspond to skew $a$-lines($a\leq r$), exceptional systems and
rulings, respectively.

As an application, we work on the monoidal transform for lines to study the
local geometry of the polytope $(r-4)_{21}$. And we show Gieser transformation
and Bertini transformation induce a symmetry of polytopes $3_{21}\ $and
$4_{21}$, respectively.

\newpage

\end{abstract}

\footnotetext{1991 Mathematics Subject Classification : 51M20,14J26,22E99}

\section{ Introduction}

The celebrated Dynkin diagrams appear as the key ingredients in many
mathematical research areas. In the geometry of polytopes, they represent the
dihedral angles between the hyperplanes generating the polytopes, and in the
algebraic geometry of surfaces, they are the intersections between the simple
roots generating a root space. In fact, the diagrams in each of above
researches correspond to the relationships presenting symmetry groups which
commonly appear in each study on the objects represented by the graphs. In
particular, the Dynkin diagrams of the Lie groups $E_{r}$ $3\leq r\leq8$
correspond to both the Weyl groups $W(S_{r})$ of del Pezzo surfaces $S_{r}$
and the symmetry groups of semiregular $E_{r}$-polytopes $(r-4)_{21}$, which
is also known as Gosset polytopes. Therefore, there are natural
correspondences between the geometry of the del Pezzo surface and the geometry
of the $(r-4)_{21}$ polytope. This article explores the correspondences
between del Pezzo surfaces and $(r-4)_{21}$ polytopes.

The del Pezzo surfaces are smooth irreducible surfaces $S_{r}$ whose
anticanonical class $-$ $K_{S_{r}}$ is ample. We can construct the del Pezzo
surfaces by blowing up $r\leq$ $8$ points from $\mathbb{P}^{2}$ unless it is
$\mathbb{P}^{1}\times\mathbb{P}^{1}$. In particular, it is very well known
that there are $27$ lines on a cubic surface $S_{6}$ and the configuration of
these lines is acted by the Weyl group $E_{6}$(\cite{Demaz}\cite{Dolgachev}%
\cite{Hartsh}). The set of $27$-lines in $S_{6}$ are bijective to the set of
vertices of a Gosset $2_{21}$ polytope, i.e. an $E_{6}$-polytope. The similar
correspondences were found for the $28$-bitangents in $S_{7}\ $and the
tritangent planes for $S_{8}$. The bijection between lines in $S_{6}$ and
vertices in $2_{21}$ was applied to study the geometry of $2_{21}$ by
Coxeter(\cite{Coxeter4}). And the complete list (see \cite{Manni}) of
bijections between the divisor classes containing lines and vertices is well
known and applied in many different research fields. These divisor classes
which are also called $lines\ $play key roles in this article.

The lines in del Pezzo surfaces are studied in many different directions.
Recently, Leung and Zhang relate the configurations of the lines to the
geometry of the line bundles over del Pezzo surfaces via the representation
theory (\cite{Leung1}\cite{Leung-Zhang}). Another interesting researches
driven from the lines in del Pezzo surfaces and their symmetry groups can be
found in \cite{Batyrev-Popov}, \cite{Friedman-Morgan} and \cite{Manv}.

A $line\ $in $Pic\ S_{r}\ $is equivalently a divisor class $l\ $with
$l^{2}=-1\ $and $K_{S_{r}}\cdot l=-1$. We observe that the Weyl group
$W(S_{r})\ $acts as an affine reflection group on the affine hyperplane given
by $D\cdot K_{S_{r}}=-1$. Furthermore, $W(S_{r})\ $acts on the set of lines in
$Pic\ S_{r}$. Therefrom, we construct a Gosset polytope $(r-4)_{21}\ $in
$Pic\ S_{r}\otimes\mathbb{Q}\ $whose vertices are exactly the lines in
$Pic\ S_{r}$. For a Gosset polytope $(r-4)_{21}$, subpolytopes are regular
simplexes excepts the facets which consist of $(r-3)$-simplexes and
$(r-3)$-crosspolytopes. Since the subpolytopes in $(r-4)_{21}$ are basically
configurations of vertices, we obtain natural characterization of subpolytopes
in $(r-4)_{21}\ $as divisor classes in $Pic\ S_{r}$.

Now, we want to use the algebraic geometry of del Pezzo surfaces to identify
the divisor classes corresponding to the subpolytopes in $(r-4)_{21}$. For
this purpose, we consider the following set of divisor classes which are
called\textit{ }$skew$\textit{ }$a$\textit{-}$lines$\textit{, }$exceptional$%
\textit{ }$systems$ and\ \textit{rulings }in $Pic\ S_{r}$.

The \textit{skew }$a$\textit{-line} is an extension of the definition of lines
in $S_{r}$. We show that each skew $a$-line represents an $(a-1)$-simplex in
an $(r-4)_{21}$ polytope. In fact, the skew $a$-lines also holds $D^{2}=-a$
and $D\cdot K_{S_{r}}=-a$. Furthermore the divisors with these conditions are
equivalently skew $a$-lines for $a\leq3$.

The \textit{exceptional system} is a divisor class in $Pic\ S_{r}$ whose
linear system gives a regular map from $S_{r}$ to $\mathbb{P}^{2}$. As this
regular map corresponds to a blowing up from $\mathbb{P}^{2}$ to $S_{r}$,
naturally we relate exceptional systems to $(r-1)$-simplexes in $(r-4)_{21}$
polytopes, which is the one of two types of facets appearing in $(r-4)_{21}$
polytopes. We show the set $\mathcal{E}_{r}$ is bijective to the set of the
$(r-1)$-simplexes in $(r-4)_{21}$ polytopes, for $3\leq r\leq$ $7$.

The \textit{ruling} is a divisor class in $Pic\ S_{r}$ which gives a fibration
of $S_{r}$ over $\mathbb{P}^{1}$. And we show that the $F_{r}$ is bijective to
the set of $(r-1)$-crosspolytopes in the $(r-4)_{21}$ polytope. Furthermore,
we explain the relationships between lines and rulings according to the
incidence between the vertices and $(r-1)$-crosspolytopes. This leads us that
a pair of proper crosspolytopes in the $(r-4)_{21}$ give the blowing down maps
from $S_{r}$ to $\mathbb{P}^{1}\times$ $\mathbb{P}^{1}$.

After proper comparison between divisor classes obtained from the geometry of
the polytope $(r-4)_{21}\ $and those given by the geometry of a del Pezzo
surface, we come by the following correspondences.%

\[%
\begin{tabular}
[c]{|l|l|}\hline
del Pezzo surface $S_{r}$ & $E$-semiregular polytopes $(r-4)_{21}%
$\\\hline\hline
lines & vertices\\\hline
skew $a$-lines $1\leq a\leq r$ & $\left(  a-1\right)  $-simplexes $1\leq a\leq
r$\\\hline
exceptional systems & $\left(  r-1\right)  $-simplexes ($r<8$)\\\hline
rulings & $\left(  r-1\right)  $-crosspolytopes\\\hline
\end{tabular}
\]

The nature of the above correspondences is macroscopic, and that we need a
microscopic explanation of the correspondences to decode the local geometry of
the $(r-4)_{21}$ polytopes. Thus, we consider the monoidal transform for lines
on del Pezzo surfaces and describe the local geometry of the $(r-4)_{21}$
polytopes. This blowing up procedure on lines can be applied to rulings to get
a useful recursive description. This will be discussed along the corresponding
geometry on the polytope $(r-4)_{21}$ in \cite{Clinger-Lee}\cite{Lee}.

As another application, we consider the pairs of lines in $Pic\ S_{7}\ $(resp.
$Pic\ S_{8}$)\ with intersection $2\ $(resp. $3$)\ which are related to the
$28$ bitangents (resp. tritangent plane). And we define Gieser transformation
(resp. Bertini transformation)\ on the polytope $3_{21}\ $(resp. $4_{21}%
$)\ and show that this is another symmetry.

The researches on regular and semiregular polytopes along the Coxeter-Dynkin
diagrams\ have a long history which may be well known only as facts. So we
begin with the preliminaries on the theories of the regular and semiregular
polytopes in the next section.

\section{\label{Section-Polytope}Regular and Semiregular Polytopes}

In this article, we deal with polytopes with highly nontrivial symmetries, and
their symmetry groups play key roles along the corresponding Coxeter-Dynkin
diagrams. In this section, we revisit the general theory of regular and
semiregular polytopes according to their symmetry groups and Coxeter-Dynkin
diagrams. Especially, we consider a family of semiregular polytopes known as
Gosset figures ($k_{21}$ according to Coxeter). The combinatorial data of
Gosset figures along the group actions will be used everywhere in this
article. For further detail of the theory, readers consult Coxeter's papers
\cite{Coxeter}\cite{Coxeter1}\cite{Coxeter2}\cite{Coxeter3}.

Let $P_{n}$\ be a convex $n$-polytope in an $n$-dimensional Euclidean space.
For each vertex $O$, the midpoints of all the edges emanating from a vertex
$O$ in $P_{n}$ form an $(n-1)$-polytope if they lie in a hyperplane. We call
this $(n-1)$-polytope the \textit{vertex figure} of $P_{n}$ at $O$.

A polytope $P_{n}$ ($n>2$) is said to \textit{regular} if its facets are
regular and there is a regular vertex figure at each vertex. When $n=2$, a
polygon $P_{2}$ is regular if it is equilateral and equiangular. Naturally,
the facets of regular $P_{n}$ are all congruent, and the vertex figures are
all the same.

We consider two classes of regular polytopes.

(1) A \textbf{regular simplex} $\alpha_{n}$ is an $n$-dimensional simplex with
equilateral edges. For example, $\alpha_{1}$ is a line-segment, $\alpha_{2}$
is an equilateral triangle, and $\alpha_{3}$ is a tetrahedron. Note
$\alpha_{n}$ is a pyramid based on $\alpha_{n-1}$. Thus the facets of a
regular simplex $\alpha_{n}$ is a regular simplex $\alpha_{n-1}$, and the
vertex figure of $\alpha_{n}$ is also $\alpha_{n-1}$. Furthermore, the
symmetry group of $\alpha_{n}$ is the Coxeter group $A_{n}$ with order
$\left(  n+1\right)  !$.

(2) A \textbf{crosspolytope} $\beta_{n}$ is an $n$-dimensional polytope whose
$2n$-vertices are the intersects between an $n$-dimensional Cartesian
coordinate frame and a sphere centered at the origin. For instance, $\beta
_{1}$ is a line-segment, $\beta_{2}$ is a square, and $\beta_{3}$ is an
octahedron. Note $\beta_{n}$ is a bipyramid based on $\beta_{n-1}$, and the
$n$-vertices in $\beta_{n}$ form $\alpha_{n-1}$ if the choice is made as one
vertex from each Cartesian coordinate line. So the vertex figure of a
crosspolytope $\beta_{n}$ is also a crosspolytope $\beta_{n-1}$, and the
facets of $\beta_{n}$ is $\alpha_{n-1}$. And the symmetry group of $\beta_{n}$
is the Coxeter group $D_{n}$with order $2^{n-1}n!$.

$\bigskip$

A polytope $P_{n}$ is called \textit{semiregular} if its facets are regular
and its vertices are equivalent, namely, the symmetry group of $P_{n}$ acts
transitively on the vertices of $P_{n}$.

Here, we consider the semiregular $k_{21}$\ polytopes discovered by
\textit{Gosset }which are $(k+4)$-dimensional polytopes whose symmetry groups
are the Coxeter group $E_{k+4}$. Note that the vertex figure of $k_{21}$ is
$\left(  k-1\right)  _{21}$ and the facets of $k_{21}$ are regular simplexes
$\alpha_{k+3}$ and crosspolytopes $\beta_{k+3}$. The list of $k_{21}$
polytopes is following.%

\[%
\begin{tabular}
[c]{|l|l|l|l|}\hline
$\ \ k$ & $E_{k+4}$ & order of $E_{k+4}$ & $k_{21}$-polytopes\\\hline\hline
$-1$ & $A_{1}\times$$A_{2}$ & $12$ & triangular prism\\\hline
$\ \ 0$ & $A_{4}$ & $5!$ & rectified 5-cell\\\hline
$\ \ 1$ & $D_{5}$ & $2^{4}5!$ & demipenteract\\\hline
$\ \ 2$ & $E_{6}$ & $72\times6!$ & $E_{6}$-polytope\\\hline
$\ \ 3$ & $E_{7}$ & $8\times9!$ & $E_{7}$-polytope\\\hline
$\ \ 4$ & $E_{8}$ & $192\times10!$ & $E_{8}$-polytope\\\hline
\end{tabular}
\]

The \textit{Coxeter groups} are reflection groups generated by the reflections
with respect to hyperplanes (called mirrors), and the \textit{Coxeter-Dynkin
diagrams} of Coxeter groups are labeled graphs where their nodes are indexed
mirrors and the labels on edges present the order $n$ of dihedral angle
$\frac{\pi}{n}$ between two mirrors. If two mirrors are perpendicular, namely
$n=2$, no edge joins two nodes presenting the mirrors because there is no
interaction between the mirrors. Since the dihedral angle $\frac{\pi}{3}$
appears very often, we only label the edges when the corresponding order
$n>3$. Each Coxeter-Dynkin diagram contains at least one ringed node which
represents an active mirror, i.e. there is a point off the mirror, and the
constructing a polytope begins with reflecting the point through the active mirror.

We call the Coxeter-Dynkin diagram of $\alpha_{n}$ (respectively $\beta_{n}$
and $k_{21}$) with the Coxeter group $A_{n}$ (respectively $D_{n}$ and $E_{n}%
$) $A_{n}$-type (respectively $D_{n}$- and $E_{n}$-type), and each
Coxeter-Dynkin diagram of $A_{n}$, $D_{n}$ and $E_{n}$-type has only one
ringed node and no labeled edges. For this case, the following simple
procedure using the Coxeter-Dynkin diagram describes possible subpolytopes and
calculates total numbers of them (see also \cite{Coxeter}\cite{Coxeter2}).

The Coxeter-Dynkin diagram of each subpolytope $P^{\prime}$ is a connected
subgraph $\Gamma$ containing the ringed node. And the subgraph obtained by
taking off all the nodes joined with the subgraph $\Gamma$ represents the
isotropy group $G_{P^{\prime}\ }$of $P^{\prime}$. Furthermore, the index
between the symmetry group $G\ $of the ambient polytope and isotropy group
$G_{P^{\prime}}$, gives the total number of such subpolytope. In particular,
by taking off the ringed node, we obtain the subgraph corresponding to the
isotropy group of a vertex, and in fact the isotropy group is the symmetry
group of the vertex figure.

\bigskip

(1) \textbf{Regular simplex }$\alpha_{n}$ with the symmetry group $A_{n}$.%

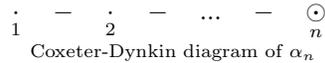
\[
\underset{\text{Coxeter-Dynkin diagram of }\alpha_{n}}{%
\begin{tabular}
[c]{lllllll}%
$\underset{1}{{\Huge \cdot}}$ & ${\Huge -}$ & $\underset{2}{{\Huge \cdot}}$ &
${\Huge -}$ & ${\Large ...}$ & ${\Huge -}$ & $\underset{n}{{\Huge \odot}}$%
\end{tabular}
\ }%
\]

The diagram of vertex figure is $A_{n-1}$-type since it is represented by the
subgraph remaining after removing the ringed node, and the facet is only
$\alpha_{n-1}$ because the subgraph of $A_{n-1}$-type is the biggest subgraph
containing the ringed node in the graph of $A_{n}$-type. Furthermore, since
all the possible subgraphs containing the ringed node are $A_{k}$-type, only
regular simplex $\alpha_{k}$, $0\leq k\leq n-1$ appears as subpolytopes. And
for each $\alpha_{k}$ in $\alpha_{n}$, the possible total number
$N_{\alpha_{k}}^{\alpha_{n}}$ is
\[
N_{\alpha_{k}}^{\alpha_{n}}=[A_{n}:A_{k}\times A_{n-k-1}]=\frac{\left(
n+1\right)  !}{(k+1)!(n-k)!}=\left(
\begin{tabular}
[c]{l}%
$n+1$\\
$k+1$%
\end{tabular}
\ \right)  .
\]

(2) \textbf{Cross polytope} $\beta_{n}$ with the symmetry group $D_{n}$.%

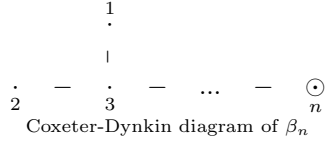
\[
\underset{\text{Coxeter-Dynkin diagram of }\beta_{n}}{%
\begin{tabular}
[c]{lllllll}
&  & $\overset{1}{{\huge \cdot}}$ &  &  &  & \\
&  & ${\Huge \shortmid}$ &  &  &  & \\
$\underset{2}{{\Huge \cdot}}$ & ${\Huge -}$ & $\underset{3}{{\Huge \cdot}}$ &
${\Huge -}$ & ${\LARGE ...}$ & ${\Huge -}$ & $\underset{n}{{\Huge \odot}}$%
\end{tabular}
}%
\]

The diagram of vertex figure is $D_{n-1}$-type because the subgraph remaining
after removing the ringed node represents $D_{n-1}$, and the facet is only
$\alpha_{n-1}$ since subgraph of $A_{n-1}$-type is the biggest subgraph
containing the ringed node in $D_{n}$-type. Only regular simplex $\alpha_{k}$,
$k=0,...,n-1$ appears as subpolytopes since the possible subgraphs containing
the ringed node are only $A_{k}$-type. And for each $\alpha_{k}\ $in
$\beta_{n}$, the possible total number $N_{\alpha_{k}}^{\beta_{n}}$ is
\[
N_{\alpha_{k}}^{\beta_{n}}=[D_{n}:A_{k}\times D_{n-k-1}]=\frac{2^{n-1}%
n!}{(k+1)!2^{n-k-2}\left(  n-k-1\right)  !}=2^{k+1}\left(
\begin{tabular}
[c]{l}%
$\ \ \ n$\\
$k+1$%
\end{tabular}
\ \ \right)  .
\]
In particular, each $\beta_{n}$ contains $N_{\alpha_{0}}^{\beta_{n}}=2n$
vertices, and these vertices form $n$-pairs with the common center.

\bigskip

(3) \textbf{Gosset polytope} $k_{21}$ with the symmetry group $E_{k+4}$,
$-1\leq k\leq4$.%

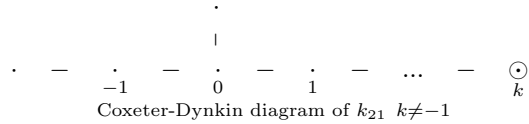
\[
\underset{\text{Coxeter-Dynkin diagram of }k_{21}~k\not =-1}{%
\begin{tabular}
[c]{lllllllllll}
&  &  &  & ${\Huge \cdot}$ &  &  &  &  &  & \\
&  &  &  & ${\Huge \shortmid}$ &  &  &  &  &  & \\
${\Huge \cdot}$ & ${\Huge -}$ & $\underset{-1}{{\Huge \cdot}}$ & ${\Huge -}$ &
$\underset{0}{{\Huge \cdot}}$ & ${\Huge -}$ & $\underset{1}{{\Huge \cdot}}$ &
${\Huge -}$ & $...$ & ${\Huge -}$ & $\underset{k}{{\Huge \odot}}$%
\end{tabular}
\ }%
\]

For $k\neq-1$, the diagram of vertex figure is $E_{k+3}$-type and the facets
are the regular simplex $\alpha_{k+3}$ and the crosspolytope $\beta_{k+3}$
since the subgraphs of $A_{k+3}$-type and $D_{k+3}$-type appear as the biggest
subgraph containing the ringed node in $E_{k+3}$-type. But all the lower
dimensional subpolytopes are regular simplexes.

Case $k=-1$ is a bit different with other cases since there are two ringed
nodes.
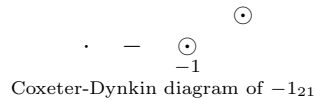
\[
\underset{\text{Coxeter-Dynkin diagram of }-1_{21}\ }{%
\begin{tabular}
[c]{llll}
&  &  & ${\Huge \odot}$\\
${\Huge \cdot}$ & ${\Huge -}$ & $\underset{-1}{{\Huge \odot}}$ &
\end{tabular}
}%
\]

The vertex figure is an isosceles triangle instead of an equilateral triangle
because the corresponding diagram is obtained by taking off a ringed node in
the $A_{2}$-type subgraph. And the facets are the regular triangle $\alpha
_{2}$ given by the $A_{2}$-type subgraph and the square $\beta_{2}$ given by
the subgraph taking off the unringed node.

As above, we can calculate the total number of subpolytopes in $k_{21}\ $by
using Coxeter-Dynkin diagram. For instance, we calculate $2_{21}$. After
removing the ringed node labelled $2$, we obtain a subgraph of $E_{5}$-type,
and therefore the vertex figure of $2_{21}\ $is $1_{21}$. Since the subgraphs
of $A_{5}$-type and $D_{5}$-type are all the possible biggest subgraphs in the
Coxeter-Dynkin diagram of $2_{21}$, there are two types of facets in $2_{21}$,
which are $5$-simplexes and $5$-crosspolytopes, respectively. And all other
subpolytopes in $2_{21}\ $are simplexes for the same reason. In the following
calculation for $2_{21}$, the nodes marked by $\ast\ \ $represent deleted
nodes.
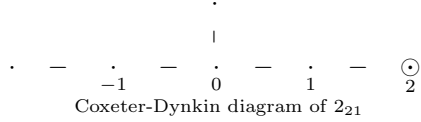
\[
\underset{\text{Coxeter-Dynkin diagram of }2_{21}}{%
\begin{tabular}
[c]{lllllllll}
&  &  &  & ${\Huge \cdot}$ &  &  &  & \\
&  &  &  & ${\Huge \shortmid}$ &  &  &  & \\
${\Huge \cdot}$ & ${\Huge -}$ & $\underset{-1}{{\Huge \cdot}}$ & ${\Huge -}$ &
$\underset{0}{{\Huge \cdot}}$ & ${\Huge -}$ & $\underset{1}{{\Huge \cdot}}$ &
${\Huge -}$ & $\underset{2}{{\Huge \odot}}$%
\end{tabular}
\ }%
\]

$\ $(a) Vertices in $2_{21}\ $: $N_{\alpha_{0}}^{2_{21}}=[E_{6}:E_{5}]=27$
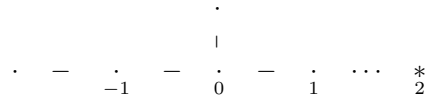
\[%
\begin{tabular}
[c]{lllllllll}
&  &  &  & ${\Huge \cdot}$ &  &  &  & \\
&  &  &  & ${\Huge \shortmid}$ &  &  &  & \\
${\Huge \cdot}$ & ${\Huge -}$ & $\underset{-1}{{\Huge \cdot}}$ & ${\Huge -}$ &
$\underset{0}{{\Huge \cdot}}$ & ${\Huge -}$ & $\underset{1}{{\Huge \cdot}}$ &
${\large \cdots}$ & $\underset{2}{\ast}$%
\end{tabular}
\]

(b) $1$-simplexes(edges) in $2_{21}$ : $N_{\alpha_{1}}^{2_{21}}=[E_{6}%
:A_{1}\times E_{4}]=216$%
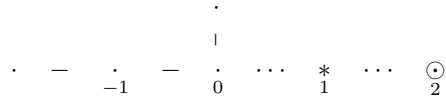
\[%
\begin{tabular}
[c]{lllllllll}
&  &  &  & ${\Huge \cdot}$ &  &  &  & \\
&  &  &  & ${\Huge \shortmid}$ &  &  &  & \\
${\Huge \cdot}$ & ${\Huge -}$ & $\underset{-1}{{\Huge \cdot}}$ & ${\Huge -}$ &
$\underset{0}{{\Huge \cdot}}$ & ${\large \cdots}$ & $\underset{1}{\ast}$ &
${\large \cdots}$ & $\underset{2}{{\Huge \odot}}$%
\end{tabular}
\ \
\]

(c) $2$-simplexes(faces) in $2_{21}\ $: $N_{\alpha_{2}}^{2_{21}}=[E_{6}%
:A_{2}\times E_{3}]=720$%
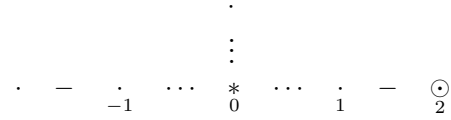
\[%
\begin{tabular}
[c]{lllllllll}
&  &  &  & ${\Huge \cdot}$ &  &  &  & \\
&  &  &  & ${\large \vdots}$ &  &  &  & \\
${\Huge \cdot}$ & ${\Huge -}$ & $\underset{-1}{{\Huge \cdot}}$ &
${\large \cdots}$ & $\underset{0}{\ast}$ & ${\large \cdots}$ & $\underset
{1}{{\Huge \cdot}}$ & ${\Huge -}$ & $\underset{2}{{\Huge \odot}}$%
\end{tabular}
\
\]

(d) $3$-simplexes(cells) in $2_{21}$ : $N_{\alpha_{3}}^{2_{21}}=[E_{6}%
:A_{3}\times A_{1}]=1080$%
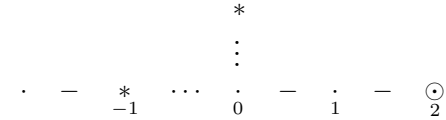
\[%
\begin{tabular}
[c]{lllllllll}
&  &  &  & $\ast$ &  &  &  & \\
&  &  &  & ${\large \vdots}$ &  &  &  & \\
${\Huge \cdot}$ & ${\Huge -}$ & $\underset{-1}{\ast}$ & ${\large \cdots}$ &
$\underset{0}{{\Huge \cdot}}$ & ${\Huge -}$ & $\underset{1}{{\Huge \cdot}}$ &
${\Huge -}$ & $\underset{2}{{\Huge \odot}}$%
\end{tabular}
\
\]

(e) $4$-simplexes in $2_{21}$ :\ $N_{\alpha_{4}}^{2_{21}}=[E_{6}:A_{4}\times
A_{1}]+[E_{6}:A_{4}]=648$%
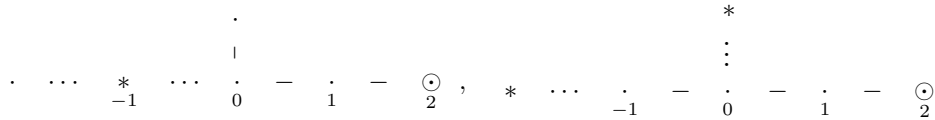
\[%
\begin{tabular}
[c]{lllllllll}
&  &  &  & ${\Huge \cdot}$ &  &  &  & \\
&  &  &  & ${\Huge \shortmid}$ &  &  &  & \\
${\Huge \cdot}$ & ${\large \cdots}$ & $\underset{-1}{\ast}$ & ${\large \cdots
}$ & $\underset{0}{{\Huge \cdot}}$ & ${\Huge -}$ & $\underset{1}{{\Huge \cdot
}}$ & ${\Huge -}$ & $\underset{2}{{\Huge \odot}}\ \ ,$%
\end{tabular}
\
\begin{tabular}
[c]{lllllllll}
&  &  &  & $\ast$ &  &  &  & \\
&  &  &  & ${\large \vdots}$ &  &  &  & \\
$\ast$ & ${\large \cdots}$ & $\underset{-1}{{\Huge \cdot}}$ & ${\Huge -}$ &
$\underset{0}{{\Huge \cdot}}$ & ${\Huge -}$ & $\underset{1}{{\Huge \cdot}}$ &
${\Huge -}$ & $\underset{2}{{\Huge \odot}}$%
\end{tabular}
\]

(f) $5$-simplexes in $2_{21}$ : $N_{\alpha_{5}}^{2_{21}}=[E_{6}:A_{5}]=72$%
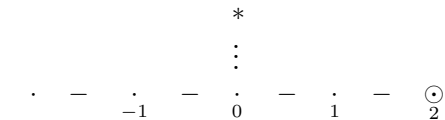
\[%
\begin{tabular}
[c]{lllllllll}
&  &  &  & $\ast$ &  &  &  & \\
&  &  &  & ${\large \vdots}$ &  &  &  & \\
${\Huge \cdot}$ & ${\Huge -}$ & $\underset{-1}{{\Huge \cdot}}$ & ${\Huge -}$ &
$\underset{0}{{\Huge \cdot}}$ & ${\Huge -}$ & $\underset{1}{{\Huge \cdot}}$ &
${\Huge -}$ & $\underset{2}{{\Huge \odot}}$%
\end{tabular}
\]

(g) $5$-crosspolytopes in $2_{21}$ : $N_{\beta_{5}}^{2_{21}}=[E_{6}%
:D_{5}]=27\ $%
\[%
\begin{tabular}
[c]{lllllllll}
&  &  &  & ${\Huge \cdot}$ &  &  &  & \\
&  &  &  & ${\Huge \shortmid}$ &  &  &  & \\
$\ast$ & ${\large \cdots}$ & $\underset{-1}{{\Huge \cdot}}$ & ${\Huge -}$ &
$\underset{0}{{\Huge \cdot}}$ & ${\Huge -}$ & $\underset{1}{{\Huge \cdot}}$ &
${\Huge -}$ & $\underset{2}{{\Huge \odot}}\ \ .$%
\end{tabular}
\]

$\ $

As we apply the same procedure to the other $E$-polytopes, we get the
following table.%

\[
\underset{\ \ \text{Numbers of subpolytopes in }k_{21}}{%
\begin{tabular}
[c]{|l||l|l|l|l|l|l|}\hline
$E_{k+4}$-polytope($k_{21}$) & $-1_{21}$ & $0_{21}$ & $1_{21}$ & $2_{21}$ &
$3_{21}$ & $4_{21}$\\\hline\hline
$\beta_{k+3}$ & $\ \ 3$ & $5$ & $10$ & $27$ & $126$ & $2160$\\\hline
vertex & $\ \ 6$ & $10$ & $16$ & $27$ & $56$ & $240$\\\hline
$\alpha_{1}$ & $\ \ 9$ & $30$ & $80$ & $216$ & $756$ & $6720$\\\hline
$\alpha_{2}$ & $\ \ 2$ & $30$ & $160$ & $720$ & $4032$ & $60480$\\\hline
$\alpha_{3}$ &  & $5$ & $120$ & $1080$ & $10080$ & $241920$\\\hline
$\alpha_{4}$ &  &  & $16$ & $648$ & $12096$ & $483840$\\\hline
$\alpha_{5}$ &  &  &  & $72$ & $6048$ & $483840$\\\hline
$\alpha_{6}$ &  &  &  &  & $576$ & $207360$\\\hline
$\alpha_{7}$ &  &  &  &  &  & $17280$\\\hline
\end{tabular}
\ \ }%
\]

\bigskip

\section{\label{SectionDelPezzo}Del Pezzo surfaces $S_{r}$}

A del Pezzo surface is a smooth irreducible surface whose anticanonical class
$-$ $K_{S}$ is ample. It is well known that a del Pezzo surface $S_{r}$,
unless it is $\mathbb{P}^{1}\times\mathbb{P}^{1}$, can be obtained from
$\mathbb{P}^{2}$ by blowing up $r\leq$ $8$ points in generic positions,
namely, no three points are on a line, no six points are on a conic, and for
$r=8$, not all of them are on a plane curve whose singular point is one of
them (see \cite{Demaz}\cite{Hartsh}\cite{Manni}).

\textbf{Notation }: We do not use different notations for the divisors and the
corresponding classes in Picard group unless there is confusion.

\bigskip

We denote such a del Pezzo surface by $S_{r}$ and the corresponding blow up by
$\pi_{r}:$ $S_{r}\rightarrow$ $\mathbb{P}^{2}$. And $K_{S_{r}}^{2}$ $=9-r$ is
called the degree of the del Pezzo surface. Each exceptional curve and the
corresponding class given by blowing up is denoted by $e_{i}$, and both the
class of $\pi_{r}^{\ast}\left(  h\right)  $ in $S_{r}$ and the class of a line
$h$ in $\mathbb{P}^{2}$ are referred as $h$. Then, we have
\[
h^{2}=1\text{, }h\cdot e_{i}=0\text{, }e_{i}\cdot e_{j}=-\delta_{ij}%
\ \ \text{for }1\leq i,j\leq r,
\]
and the Picard group\ of $S_{r}$ is
\[
Pic\ S_{r}\simeq\mathbb{Z}h\oplus\mathbb{Z}e_{1}\oplus...\oplus\mathbb{Z}e_{r}%
\]
with the signature $(1,-r)$. And $K_{S_{r}}=-3h+\sum_{i=1}^{r}e_{i}$.

For any irreducible curve $C$ on a Del Pezzo surface $S_{r}$, we have $C\cdot
K_{S_{r}}$ $<0$ since $-K_{S_{r}}$ is ample. Furthermore, if the curve $C$ has
a negative self-intersection, $C$ must be a smooth rational curve with
$C^{2}=-1$ by the adjunction formula.

The ample\ $-K_{S_{r}}\ $on a del Pezzo surface $S_{r}\ $is very beneficial to
deal with $Pic\ S_{r}$. The inner product given by the intersection on
$Pic\ S_{r}\ $induces a negative definite metric on $\left(  \mathbb{Z}%
K_{S_{r}}\right)  ^{\perp}$ in $Pic\ S_{r}\ $where we can also define natural reflections.

To define reflections on $\left(  \mathbb{Z}K_{S_{r}}\right)  ^{\perp}$ in
$Pic\ S_{r}$, we consider a root system%

\[
R_{r}:=\{d\in Pic\ S_{r}\mid d^{2}=-2,\ d\cdot K_{S_{r}}=0\}\text{,}%
\]
with simple roots
\[
d_{0}=h-e_{1}-e_{2}-e_{3},d_{i}=e_{i}-e_{i+1},\ 1\leq i\leq r-1.
\]
Each element $d$ in $R_{r}$ defines a reflection on $\left(  \mathbb{Z}%
K_{S_{r}}\right)  ^{\perp}$ in $Pic\ S_{r}$
\[
\sigma_{d}(D):=D+\left(  D\cdot d\right)  d\ \ \text{for\ }D\in\left(
\mathbb{Z}K_{S_{r}}\right)  ^{\perp}%
\]
and the corresponding Weyl group $W(S_{r})$ is $E_{r}$ where $3\leq r\leq8$
with the Dynkin diagram%
\[
\underset{\text{Dynkin diagram of }E_{r}\ r\geq3}{%
\begin{tabular}
[c]{lllllllllll}
&  &  &  & $\overset{d_{0}}{{\Huge \cdot}}$ &  &  &  &  &  & \\
&  &  &  & $\ {\Huge \shortmid}$ &  &  &  &  &  & \\
$\underset{d_{1}}{{\Huge \cdot}}$ & ${\Huge -}$ & $\underset{d_{2}%
}{{\Huge \cdot}}$ & ${\Huge -}$ & $\underset{d_{4}}{{\Huge \cdot}}$ &
${\Huge -}$ & $\underset{d_{5}}{{\Huge \cdot}}$ & ${\Huge -}$ & ${\Large ...}$
& ${\Huge -}$ & $\underset{d_{r-1}}{{\Huge \cdot}}$%
\end{tabular}
\ \ \ \ \ \ }%
\]

The definition of the reflection $\sigma_{d}$ on $\left(  \mathbb{Z}K_{S_{r}%
}\right)  ^{\perp}$ can be used to obtain a transformation both on
$Pic\ S_{r}$ and on $Pic\ S_{r}\otimes\mathbb{Q\simeq Q}h\oplus\mathbb{Q}%
e_{1}\oplus...\oplus\mathbb{Q}e_{r}$ via the linear extension of the
intersections of divisors in $Pic\ S_{r}$. Here $Pic\ S_{r}\otimes\mathbb{Q}$
is a vector space with the signature $(1,-r)$.

\bigskip

\textbf{Affine hyperplanes and the reflection groups}

Later on, we deal with divisor classes $D\ $satisfying $D\cdot K_{S_{r}%
}=\alpha$,\ $D^{2}=\beta\ $where $\alpha\ $and $\beta\ $are integers along the
action of Weyl group $W(S_{r})$. Here, we know $W(S_{r})\ $is generated by the
reflections on $\left(  \mathbb{Z}K_{S_{r}}\right)  ^{\perp}\ $given by simple
roots. To extend the action of $W(S_{r})$ properly, we want to show that these
reflections are defined on $Pic\ S_{r}\ $and preserve the above equations.
Furthermore, we see that $W(S_{r})\ $acts as a reflection group on the set of
divisor classes with $D\cdot K_{S_{r}}=\alpha$.

We consider an \textit{affine hyperplane section} in $Pic\ S_{r}%
\otimes\mathbb{Q}$ defined by
\[
\tilde{H}_{b}:=\{D\in Pic\ S_{r}\ \otimes\mathbb{Q}\mid-D\cdot K_{S_{r}}=b\}
\]
where $b$ is an arbitrary real number and an affine hyperplane section
$H_{b}:=\tilde{H}_{b}\cap Pic\ S_{r}\ $in $Pic\ S_{r}$. Since $-K_{S_{r}}\ $is
ample, we are interested in $b\geq0$.

By the fact that $K_{S_{r}}^{2}=9-r>0$, $3\leq r\leq$ $8$ and Hodge index
theorem%
\[
0=(K_{S_{r}}\cdot(D_{1}-D_{2}))^{2}\geq K_{S_{r}}^{2}(D_{1}-D_{2}%
)^{2}\ \text{, }D_{1},D_{2}\in H_{b}\text{,}%
\]
the inner product on $Pic\ S_{r}$ induces a negative definite metric on
$H_{b}$. As a matter of fact, the induced metric is defined on $Pic\ S_{r}%
\otimes\mathbb{Q}$, and we can also consider the induced norm by fixing a
center $\frac{b}{9-r}K_{S_{r}}$ in the affine hyperplane section $-D\cdot
K_{S_{r}}=b$ in $Pic\ S_{r}\otimes\mathbb{Q}$. This norm is also negative definite.

\begin{lemma}
\label{LMAHyperp-Refl}(1)Let $\tilde{H}_{b}\ (b\geq0)$ be an affine hyperplane
section in $Pic\ S_{r}\otimes\mathbb{Q}$ defined above\ and $\frac{b}%
{9-r}K_{S_{r}}$ be a center on the affine hyperplane section. The classes $D$
in $H_{b}=\tilde{H}_{b}\cap$ $Pic\ S_{r}\ $with a fixed self-intersection are
on a sphere with the center $\frac{b}{9-r}K_{S_{r}}$ in $\tilde{H}_{b}$.

(2) For each root $d$ in $R_{r}$, the corresponding reflection $\sigma_{d}$
defined on $Pic\ S_{r}\otimes\mathbb{Q}$ is an isometry preserving $K_{S_{r}}$
and acts as a reflection on each hyperplane section $\tilde{H}_{b}\ $with the
center $\frac{b}{9-r}K_{S_{r}}$.
\end{lemma}

\textbf{Proof}: (1) Consider
\[
\left(  D-\frac{b}{9-r}K_{S_{r}}\right)  ^{2}=D^{2}-\frac{2b}{9-r}D\cdot
K_{S_{r}}+\frac{b^{2}}{\left(  9-r\right)  ^{2}}K_{S_{r}}^{2}=D^{2}%
-\frac{b^{2}}{\left(  9-r\right)  }\leq0
\]
and the last inequality is given by the Hodge index theorem
\[
b^{2}=\left(  D\cdot K_{S_{r}}\right)  ^{2}\geq D^{2}K_{S_{r}}^{2}%
=D^{2}\left(  9-r\right)  .
\]

(2) Each root $d\ $in $R_{r}\ $satisfies $d\cdot$ $K_{S_{r}}=0\ $and
$d^{2}=-2$. Therefore, we have
\[
\sigma_{d}(K_{S_{r}})=K_{S_{r}}+(d\cdot K_{S_{r}})d=K_{S_{r}}\text{,}%
\]
and for each $D_{1},D_{2}\ \in Pic\ S_{r}\otimes\mathbb{Q}$
\[
\sigma_{d}(D_{1})\cdot\sigma_{d}(D_{2})=\left(  D_{1}+(d\cdot D_{1})d\right)
\cdot\left(  D_{2}+(d\cdot D_{2})d\right)  =D_{1}\cdot D_{2}\text{.}%
\]
Furthermore, for each class $D$ in $Pic\ S_{r}\otimes\mathbb{Q}$, the
self-intersection $D^{2}$ and $D\cdot K_{S_{r}}$ are invariant under
$\sigma_{d}$. This implies $\sigma_{d}\ $acts on the hyperplane section
$\tilde{H}_{b}$. Moreover, the hyperplane in $Pic\ S_{r}\otimes\mathbb{Q\ }%
$preserved by the action of $\sigma_{d}\ $is given by an equation $d\cdot
D=0\ $for $D\ \in Pic\ S_{r}\otimes\mathbb{Q}$, and the each center $\frac
{b}{9-r}K_{S_{r}}\ $of $\tilde{H}_{b}\ $is in this hyperplane. And because
each class $D$ in $\tilde{H}_{b}\ $can be written as
\[
D=D_{3}+\frac{b}{9-r}K_{S_{r}}\ \ \text{for some\ }D_{3}\in\tilde{H}%
_{0}\text{,}%
\]
we have%
\[
\sigma_{d}(D)=\sigma_{d}(D_{3}+\frac{b}{9-r}K_{S_{r}})=\sigma_{d}(D_{3}%
)+\frac{b}{9-r}K_{S_{r}}\in\tilde{H}_{0}+\frac{b}{9-r}K_{S_{r}}=\tilde{H}%
_{b}\text{.}%
\]
Since $\sigma_{d}\ $is a reflection on $\tilde{H}_{0}$, we can derive a fact
that the isometry $\sigma_{d}\ $acts as an affine reflection on $\tilde{H}%
_{b}\ $for the center $\frac{b}{9-r}K_{S_{r}}$. $\ \ \ {\LARGE \blacksquare}$

\bigskip

Generically, the hyperplanes in $Pic\ S_{r}\otimes\mathbb{Q\ }$induce affine
hyperplanes in $\tilde{H}_{b}\ $and they may not share a common point. But the
reflection hyperplane of each reflection $\sigma_{d}\ $in $Pic\ S_{r}%
\otimes\mathbb{Q\ }$gives a hyperplane in $\tilde{H}_{b}\ $containing the
center because it is given by a condition $K_{S_{r}}\cdot d=0$. Therefrom, the
above lemma gives the following corollary.

\begin{corollary}
\label{Coro-Affine-Relf}The affine reflections on $\tilde{H}_{b}$ given by
simple roots in $R_{d}\ $generate the Weyl group $W(S_{r})$.
\end{corollary}

\bigskip

\textbf{Remark} : The Weyl group $W(S_{r})$, generated by the simple roots in
$R_{d}$ is also preserves the self-intersection of each divisor in
$Pic\ S_{r}$ and acts on each $H_{b}=Pic\ S_{r}\cap$ $\tilde{H}_{b}\ $as an
affine reflection group. According to \cite{Friedman-Morgan}, the Weyl group
$W(S_{r})$ is the isotropy group of $K_{S_{r}}$ in the automorphism group of
$Pic\ S_{r}$.

\bigskip

\section{\label{Sec-GosinPic}Gosset Polytopes $(r-4)_{21}$ in $Pic\ S_{r}%
\otimes\mathbb{Q}$}

In this section, we identify a special classes in $Pic\ S_{r}$, which is known
as a line, and$\ $construct Gosset polytopes $(r-4)_{21}$ in $Pic\ S_{r}%
\otimes\mathbb{Q}\ $as the convex hull of the set of lines. And we study the
divisor classes representing subpolytopes in $(r-4)_{21}$.

\subsection{\label{SubSec-Constuction}Lines in del Pezzo surface $S_{r}$}

The configuration of lines on a del Pezzo surface $S_{r}\ $has been driven big
attention because of its high degree of symmetry related to the Weyl group
$W(S_{r})$ of $E_{r}$-type. When $r\leq6$, the anticanonical class $-K_{S_{r}%
}\ $on del Pezzo surface $S_{r}\ $is very ample and its linear system gives an
imbedding to$\ \mathbb{P}^{9-r}\ $where $K_{S_{r}}^{2}=9-r$. And a smooth
rational curve $C$ in $S_{r}\ $is mapped to a line in $\mathbb{P}^{9-r}\ $\ if
and only if it is exceptional curve in $S_{r}$. Furthermore, the divisor class
$D\ $containing the curve $C\ $satisfies $D\cdot K_{S_{r}}=-1=D^{2}\ $and the
vice versa \cite{Manni}. Since the last numerical equivalence is true for each
of del Pezzo surfaces, we also call the divisors with the these conditions in
del Pezzo surfaces \textit{lines}. As the symmetry group of lines in the cubic
is the Weyl group $E_{6}$, the symmetry group of lines in $S_{r}$ is the Weyl
group $E_{r}$.

We define the set of lines on $Pic\ S_{r}$\ as%
\[
L_{r}:=\{l\in Pic(S_{r})\mid l^{2}=l\cdot K_{S_{r}}=-1\}.
\]
By the adjunction formula a divisor in this class represents a rational smooth
curve in $S_{r}$. By going through simple calculation, we can obtain the
number of lines in $L_{r}$, and moreover the numbers of sets of lines are
parallel to the numbers of vertices in Gosset polytopes $(r-4)_{21}$.%
\[%
\begin{tabular}
[c]{|l|l|l|l|l|l|l|}\hline
del Pezzo Surfaces & \ $\ S_{3}$ & $\ S_{4}$ & $\ S_{5}$ & $\ S_{6}$ & $S_{7}$
& $S_{8}$\\\hline
number of Lines & $\ \ 6$ & $10$ & $16$ & $27$ & $56$ & $240$\\\hline\hline
Gosset Polytopes $(r-4)_{21}$ & $-1_{21}$ & $\ 0_{21}$ & $1_{21}$ & $2_{21}$ &
$3_{21}$ & $4_{21}$\\\hline
number of Vertices & $\ \ 6$ & $10$ & $16$ & $27$ & $56$ & $240$\\\hline
\end{tabular}
\ \ \
\]

In fact, this bijection between lines and vertices is a well-known fact
(\cite{Manni}\cite{Dolgachev}). In this article, this fact induces significant
implications after our construction of Gosset polytopes in $Pic\ S_{r}\ $where
each vertex automatically represents a line.

First of all, we need to consider intersections between lines and roots in
$Pic$ $S_{r}$.

The possible intersections of the lines in $Pic$ $S_{r}$ can be obtained by
the Hodge index theorem,%

\[
\left(  K_{S_{r}}\cdot\left(  l_{1}\pm l_{2}\right)  \right)  ^{2}\geq
K_{S_{r}}^{2}\left(  l_{1}\pm l_{2}\right)  ^{2}.
\]
And we have
\[
\frac{2}{9-r}+1\geq l_{1}\cdot l_{2}\geq-1.
\]
Therefore, two distinct lines $l_{1}$ and $l_{2}$ in $Pic\ S_{r}$ can have
intersections such as
\[
l_{1}\cdot l_{2}=\left\{
\begin{tabular}
[c]{l}%
$0,1\ \ \ \ \ \ $\\
$0,1,2\ \ \ \ \ \ $\\
$0,1,2,3$%
\end{tabular}%
\begin{tabular}
[c]{l}%
$3\leq r\leq6$\\
$\ r=7$\\
$\ r=8$%
\end{tabular}
\ \right.  .
\]
Furthermore, $l_{1}\cdot l_{2}=2$ for $r=7$ and $l_{1}\cdot l_{2}=3$ for $r=8$
hold the equalities in the Hodge index theorem, and we have equivalences%
\begin{align*}
l_{1}\cdot l_{2}  &  =2\Longleftrightarrow l_{1}+l_{2}=-K_{S_{7}}\ \text{for
}r=7,\\
l_{1}\cdot l_{2}  &  =3\Longleftrightarrow l_{1}+l_{2}=-2K_{S_{8}}~\text{for
}r=8\text{.}%
\end{align*}

Recall that for a reflection $\sigma_{d}$ given by a root $d$, if $l$ is a
line $\sigma_{d}(l)$ is also a line by lemma \ref{LMAHyperp-Refl}. Moreover we
have
\[
\sigma_{d}(l)\cdot l=\left(  l+\left(  l\cdot d\right)  d\right)  \cdot
l=-1+\left(  l\cdot d\right)  ^{2}.
\]

From the above possible numbers of the intersections of lines, the possible
intersections between a line $l$ and a root $d$ are given as
\[
l\cdot d=\left\{
\begin{tabular}
[c]{l}%
$0,\pm1\ \ \ \ \ \ \ $\\
$0,\pm1,\pm2\ $%
\end{tabular}%
\begin{tabular}
[c]{l}%
$3\leq r\leq7$\\
$r=8$%
\end{tabular}
\ \ \right.  .
\]
Since $\sigma_{d}(l)\cdot l+1$ must be a square of these integers, it is easy
to see that any two lines $l_{1}$and $l_{2}$ with $l_{1}\cdot l_{2}=1$ or $2$
cannot be mapped to each other by a reflection $\sigma_{d}$ given by a root
$d$.

\begin{lemma}
\label{Lemma-reflection-Root-Lines}(1) For each line $l$ in $S_{r}$, the
reflection $\sigma_{d}$ given by a root $d$ preserves the line $l$ if and only
if $d\cdot l=0$.

(2) Any two distinct lines in $Pic\ S_{r\leq7}$ are skew if and only if there
is a reflection $\sigma_{d}$ given by a root $d$ which reflects these lines to
each other. For $S_{8}$, this statement is true with an extra condition that a
root $d$ is chosen to have intersection $1$ with one of the lines.
\end{lemma}

\bigskip Note: We call distinct lines $l_{1}$ and $l_{2}$ are \textit{skew} if
$l_{1}\cdot l_{2}=0$.

\textbf{Proof}: (1) Trivial from
\[
l=\sigma_{d}\left(  l\right)  =l+\left(  d\cdot l\right)  d.
\]

(2) If $l_{1}$and $l_{2}$ are skew, namely $l_{1}\cdot l_{2}=0$, then
$l_{1}-l_{2}$ is a root and the corresponding reflection $\sigma_{l_{1}-l_{2}%
}$ satisfies $\sigma_{l_{1}-l_{2}}\left(  l_{1}\right)  =l_{2}$ and
$\sigma_{l_{1}-l_{2}}\left(  l_{2}\right)  =l_{1}$. Conversely, when $3\leq
r\leq7$ if two distinct lines $l_{1}$and $l_{2}$ satisfies
\[
l_{2}=\sigma_{d}\left(  l_{1}\right)  =l_{1}+\left(  d\cdot l_{1}\right)
d,\ \text{for a root }d\ \text{with }d\cdot l_{1}\not =0\text{,}%
\]
then
\[
l_{1}\cdot l_{2}=-1+\left(  d\cdot l_{1}\right)  ^{2}=0
\]
according to above list of possible intersections between a line and a root.
The case $r=8\ $is similar to the other once we add the condition that $d\cdot
l_{1}=\pm1$.$\ {\LARGE \blacksquare}$

\bigskip

By lemmar\ \ref{LMAHyperp-Refl}, the action of Weyl group $W(S_{r})$ preserves
the conditions $l^{2}=l\cdot K_{S_{r}}=-1$, and therefore $W(S_{r})$ acts on
the set of lines $L_{r}$ on $S_{r}$. Furthermore, by following theorem, there
is only one orbit of $W(S_{r})$ in the set of lines $L_{r}$, and it implies
that the bijection between the set of lines and the set of vertices in the
above is more than the correspondence between sets.

\bigskip

\begin{theorem}
\label{Correspondence line and vetices}For each del Pezzo surface $S_{r}$, the
set of lines $L_{r}$ on $S_{r}$ is the set of vertices of a Gosset polytope
$(r-4)_{21}\ $in a hyperplane section $\tilde{H}_{1}$.
\end{theorem}

\textbf{Proof} : Recall the center of $\tilde{H}_{1}\ $is $\frac{-K_{S_{r}}%
}{9-r}\ $and the distance between a line and the center in $\tilde{H}_{1}$ is
$-1-\frac{1}{9-r}$. Therefore, the set of lines $L_{r}\ $sits in a sphere in
$\tilde{H}_{1}\ $with the center $\frac{-K_{S_{r}}}{9-r}$. Furthermore, the
convex hull of $L_{r}\ $in $\tilde{H}_{1}\ $is a convex polytope.\ We want to
show this polytope is a Gosset\ polytope $(r-4)_{21}$. We construct a Gosset
polytope $(r-4)_{21}\ $in the convex hull of $L_{r}$, and we show the convex
hull is same with the polytope $(r-4)_{21}$. By lemma \ref{LMAHyperp-Refl},
the set $L_{r}\ $is acted by the Weyl group $W(S_{r})$. We choose a line
$e_{r}\ $in $L_{r}\ $and consider the generators of $W(S)\ $given by simple
roots $d_{i}\ $($0\leq i\leq r-1\ $). Since $d_{i}\cdot e_{r}=0\ $except
$i=r-1$,\ the reflection given by $d_{r-1}\ $is only active among the
generators. The line $e_{r}\ $and the generators of $W(S_{r})$ give the
Coxeter-Dynkin diagram of $E_{r}$-type\ with a ringed node at $d_{r-1}$.
Therefore, via the action of $W(S_{r})\ $on $e_{r}\ $as in the section
\ref{Section-Polytope},\ we obtain a Gosset polytope $(r-4)_{21}\ $containing
$e_{r}\ $in the convex hull of $L_{r}$. Now as we know, the number of vertices
in $(r-4)_{21}\ $is the same with the number of lines in $L_{r}$, therefore
the polytope $(r-4)_{21}\ $and the convex hull of $L_{r}\ $are the same. This
concludes the theorem.\ \ ${\LARGE \blacksquare}$

\bigskip

This theorem implies that the Weyl group $W(S_{r})$ acts transitively on
$L_{r}$. And we also see that the integral classes representing lines in
$Pic(S_{r})\ $are honest vertices of $(r-4)_{21}\ $in an affine hyperplane
$H_{1}$. But we also denote a line $l\ $in $Pic\ S_{r}$ by $V_{l}\ $if it is
considered as a vertex of $(r-4)_{21}$.

\bigskip

\begin{corollary}
\ The number of lines in the del Pezzo surface $S_{r}$ is the same with the
number of vertices of the $E_{r}$-semiregular polytope. Furthermore, the Weyl
group $W(S_{r})$ acts transitively on the set of lines $L_{r}$ on $S_{r}$.
\end{corollary}

\bigskip

\textbf{Remark}: From the section \ref{Section-Polytope}, we know that the
isotropy group of a vertex of $(r-4)_{21}\ $is the same with the symmetry
group of vertex figure which is $E_{r-1}$-type. We can check this for the
Gosset polytope $(r-4)_{21}\ $in $Pic\ S_{r}\otimes\mathbb{Q}$.\ By above
corollary, we can choose an exceptional class $e_{r}\ $without losing
generality. The generators of the isotropy group of $e_{r}$ are the simple
roots in the Dynkin diagram of $E_{r}$ which are perpendicular to $e_{r}$. And
the relationships of these simple roots are presented as a subdiagram of the
Dynkin diagram of $E_{r}$ by takin off the node $d_{r-1}$. Therefore, the
isotropy group of $e_{r}$ in $W(S_{r})$ is the Weyl group of $W(S_{r-1})\ $of
$E_{r-1}$-type, and similarly the isotropy group of each line in $S_{r}$ is
conjugate to $W(S_{r-1})$.

\bigskip

\textbf{Intersections of lines and configuration of vertices}

As the configuration of lines is our main issue, the intersection between
lines characterize how the corresponding vertices are related the polytope
$(r-4)_{21}$. Here we discuss the relationship via the properties of the lines
with fixed intersections, and the complete and uniform description will be
given in the last section.

For a fixed line $l$, the isotropy group of a line $l\ $is $E_{r-1}$-type, and
more it is same with the symmetry group of the vertex figure of $(r-4)_{21}$.
In fact, the vertex set of vertex figure corresponds to the set of lines
intersecting zero to the line $l$ by the following lemma.

\begin{lemma}
\label{Edged}For each distinct lines $l_{1}$ and $l_{2}$ in $Pic\ S_{r}$, the
vertices $V_{l_{1}}$ and $V_{l_{2}}$ in the Gosset polytope $(r-4)_{21}$ in
$\tilde{H}_{1}\ $are joined by an edge if and only of $l_{1}\cdot l_{2}=0$.
\end{lemma}

\textbf{Proof}: The semiregular polytope $(r-4)_{21}\ $in $\tilde{H}_{1}$ is
the convex hull of its vertices. Therefore, a fixed vertex and the vertices
edged to it have the minimal distance among the distance between vertices.
Since the metric on the hyperplane $\tilde{H}_{1}$ is negative definite(see
section \ref{SectionDelPezzo}), the distance $\left(  l_{1}-l_{2}\right)
^{2}$ among lines in $L_{r}\ $is maximal if and only if $l_{1}\cdot l_{2}=0$.
Therefore we have the lemma. ${\LARGE \blacksquare}$

\bigskip

\textbf{Remark}: In fact, this lemma is one case of skew $2$-lines of theorem
\ref{alinevssimplex} in the next subsection which is proven by different argument.

\bigskip

For a fixed line $l$ in a del Pezzo surface $S_{r}$, $3\leq r\leq8$, we
consider a set $N_{k}(l,S_{r})$ defined as%
\[
N_{k}(l,S_{r}):=\left\{  \ l^{\prime}\ \in L_{r}\mid l^{\prime}\cdot
l=k\right\}  .
\]
\ From above lemma, we have $\left\vert N_{0}(l,S_{r})\right\vert =\left\vert
L_{r-1}\right\vert =\left[  E_{r-1}:E_{r-2}\right]  \ $where the last equality
is given by theorem \ref{Correspondence line and vetices}.

The simple comparison according to the above list of intersections of lines
leads us the following useful lemma.

\begin{lemma}
\label{lemma-int2}For each line $l\ $in $Pic\ S_{r}$ ($4\leq r\leq
8)$,\ $\left\vert N_{1}(l,S_{r})\right\vert \ $equals the number of
$(r-2)$-crosspolytopes in the polytope $(r-3)_{21}\ $, i.e. $\left[
E_{r-1}:D_{r-2}\right]  $.
\end{lemma}

\textbf{Proof}:\ Since Weyl group $W(S_{r})\ $acts transitively, we can choose
an exceptional class $e_{r}\ $and the results on this line also hold for the
other lines. When $\ 4\leq r\leq6$, a line $e_{r}\ $intersects the other lines
with $0\ $or $1$. Since $N_{0}(e_{r},S_{r})\ $is the number of lines in the
vertex figure, we have
\begin{align*}
N_{1}(e_{r},S_{r})  &  =\left\vert L_{r}\right\vert -\left\vert N_{0}%
(e_{r},S_{r})\right\vert -\left\vert \{e_{r}\}\right\vert \\
&  =\left\vert L_{r}\right\vert -\left\vert L_{r-1}\right\vert -1\ \text{.}%
\end{align*}
So $N_{1}(e_{4},S_{4})=3$,\ $N_{1}(e_{5},S_{5})=5\ $and $N_{1}(e_{6}%
,S_{6})=10$. These are exactly the numbers of crosspolytopes in $-1_{21}%
$,\ $0_{21}\ $and $1_{21}\ $respectively.

For $r=7,\ $a line $e_{7}\ $meets the other lines with intersection $0$,\ $1$
and $2$. As we saw above, $-K_{S_{7}}-e_{7}$ is only one line with
intersecting $e_{7}\ $by$\ 2$. So we have%
\begin{align*}
N_{1}(e_{7},S_{7})  &  =\left\vert L_{7}\right\vert -\left\vert N_{0}%
(e_{7},S_{7})\right\vert -\left\vert N_{2}(e_{7},S_{7})\right\vert -\left\vert
\{e_{7}\}\right\vert \\
&  =\left\vert L_{7}\right\vert -\left\vert L_{6}\right\vert -2=27\ \text{,}%
\end{align*}
and this is the number of $6$-crosspolytopes in $2_{21}$.

For $r=8,\ $a line $e_{8}\ $may have intersection $0$,\ $1$,\ $2\ $and
$3\ $with other lines.\ Later in this article, we use a transformation between
lines in $Pic\ S_{8}\ $defined by $-\left(  2K_{S_{8}}+l\right)  $, for each
line $l\ $in $L_{8}$. Observe that a line $l\ $intersects $e_{8}\ $by $0\ $if
and only if $-\left(  2K_{S_{8}}+l\right)  \ $intersects $e_{8}\ $by $2$.
Therefore $N_{0}(e_{8},S_{8})=N_{2}(e_{8},S_{8})$. Since $-2K_{S_{8}}-e_{8}$
is only one line in $L_{8}\ $intersecting $e_{8}\ $by $3$, we get
\begin{align*}
N_{1}(e_{8},S_{8})  &  =\left\vert L_{8}\right\vert -\left\vert N_{0}%
(e_{8},S_{8})\right\vert -\left\vert N_{2}(e_{8},S_{8})\right\vert -\left\vert
N_{3}(e_{8},S_{8})\right\vert -\left\vert \{e_{8}\}\right\vert \\
&  =\left\vert L_{8}\right\vert -2\left\vert L_{7}\right\vert -2=126\ \text{,}%
\end{align*}
and this is the number of $7$-crosspolytopes in $3_{21}$.
${\LARGE \blacksquare}$

\bigskip

\textbf{Remark}:\ The above results on $N_{k}(l,S_{r})\ $appear again when we
study the monoidal transforms to lines in del Pezzo surfaces.

\subsection{\label{Sebsec-geomPoly}Subpolytopes in Gosset Polytope
$(r-4)_{21}$}

By the theorem \ref{Correspondence line and vetices}, we construct the Gosset
polytope $(r-4)_{21}$ in $Pic\ S_{r}\otimes\mathbb{Q}\ $and we want to
characterize each subpolytopes in $(r-4)_{21}\ $as a class in $Pic\ S_{r}$.
There are $k$-simplexes\ $0\leq k\leq r-1$ and $(r-1)$-crosspolytopes in
$(r-4)_{21}$.

To identify each subpolytope in $(r-4)_{21}$, we want to use the barycenter of
the subpolytope. By the theorem \ref{Correspondence line and vetices}, each
vertex of the polytope $(r-4)_{21}$ represents a line in $S_{r}$, and the
honest centers of simplexes (resp. crosspolytopes ) are written as
$(l_{1}+...+l_{k})/k$ (resp.$(l_{1}^{\prime}+l_{2}^{\prime})/2$)\ in
$\tilde{H}_{1}$ which may not be elements in $Pic$ $S_{r}$. Therefore,
alternatively, we choose $(l_{1}+...+l_{k})$ as the center of a subpolytope so
that $(l_{1}+...+l_{k})$ is in $Pic$ $S_{r}$.

\bigskip

\textbf{Simplexes in }$(r-4)_{21}$

Each $a$-simplexes in $(r-4)_{21}\ $shares its edges with $(r-4)_{21}$ and
consist of $a+1\ $vertices in $(r-4)_{21}\ $joined with edges to each other.
In our way to understand, these vertices correspond to lines in $L_{r}\ $where
they are skew to each other, namely disjoint, by lemma \ref{Edged}. And we
consider the set of center of $a$-simplexes $0\leq a\leq r-1\ $in
$(r-4)_{21}\ $defined by%
\[
A_{r}^{a}:=\{D\in Pic(S_{r})\mid\text{ }D=l_{1}+...+l_{a+1}\text{,\ }%
l_{i}\ \text{disjoint lines in }L_{r}\}\text{.}%
\]

Thanks to convexity, we expect each center in $A_{r}^{a}\ $represents an
$a$-simplex in $(r-4)_{21}$, and we prove this by algebraically as follows.
Therefrom, we have a bijection between $A_{r}^{a}\ $and the set of
$a$-simplexes in $(r-4)_{21}$.

\begin{lemma}
All the center of $a$-simplexes in $(r-4)_{21}$ are distinct.
\end{lemma}

\textbf{Proof}. Since $a=0$ is trivial, we assume\ $1\leq a\leq r-1$.\ Let
$l_{i},$\ $1$\ $\leq i\leq a+1$ be skew lines of an $a$-simplex $P$ in
$(r-4)_{21}$, and assume there is another set of disjoint lines $l_{i}%
^{\prime}$\ $1\leq i\leq a+1$ of an $a$-simplex\ $P^{\prime}$ whose center is
same with the center of $P$. We observe
\[
\left(  l_{1}+...+l_{a+1}\right)  \cdot l_{a+1}^{\prime}=\left(  l_{1}%
^{\prime}+...+l_{a+1}^{\prime}\right)  \cdot l_{a+1}^{\prime}=-1.
\]
Therefore, $l_{a+1}^{\prime}$ is one of the lines $l_{i}$ $1\leq i\leq a+1$.
And by induction, the choice of lines $l_{a}$ and the choice of lines
$l_{a}^{\prime}$ are the same. Therefore, the $a$-simplexes $P\ $and
$P^{\prime}\ $are the same. This gives lemma. $\ \ \ {\LARGE \blacksquare}$

\bigskip

\textbf{Crosspolytopes in }$(r-4)_{21}$

An $(r-1)$-crosspolytope consists of $(r-1)\ $pairs of vertices where all the
pairs of vertices share a common center which is also the center of the
polytope. Furthermore, a vertex in the polytope joined to the other vertices
by edges in the polytope except a vertex making the center with the given
vertex. According to the theorem\ \ref{Correspondence line and vetices}, an
$(r-1)$-crosspolytope in $Pic\ S_{r}$ is determined by $(r-1)\ $pairs of
lines. Since these pairs shares a common center, we consider two pairs of
lines from the polytopes such as $(l_{1},l_{2})\ $and $(l_{1}^{\prime}%
,l_{2}^{\prime})\ $with $l_{1}+l_{2}=l_{1}^{\prime}+l_{2}^{\prime}\ $. Because
$l_{1}\ $is joined to $l_{1}^{\prime}\ $and $l_{2}^{\prime}$, by lemma
\ref{Edged} we have $l_{1}\cdot l_{2}=l_{1}\cdot(l_{1}^{\prime}+l_{2}^{\prime
}-l_{2})=1$. Therefore, each pair of lines in the crosspolytope has
intersection $1$.\ In the following theorem, it turns out that a couple of
lines with intersection $1\ $characterize an $(r-1)$-crosspolytope in
$(r-4)_{21}\ $. It is useful to note that one can get $l_{i}\cdot
l_{1}^{\prime}=0=l_{i}\cdot l_{2}^{\prime}\ \ i=1,2\ \ $for above pairs, and
it implies any two pairs of lines with intersection $1$ sharing a common
center correspond to two diagonal pairs of vertices in a square. Therefore, a
center cannot be shared by more than one $(r-1)$-crosspolytope.

\begin{theorem}
\label{Thm-CenterCrosspoly}For a del Pezzo surface $S_{r}$,
\[
B_{r}:=\{D\in Pic(S_{r})\mid D=l_{1}+l_{2},\ \text{where}\ \ l_{1}\ \text{,
}l_{2}\text{ lines with\ }l_{1}\cdot l_{2}=1\},
\]
is the set of the center of $(r-1)$-crosspolytopes in $(r-4)_{21}$.
\end{theorem}

\textbf{Proof:\ }From the above, we know each center of an $(r-1)$%
-crosspolytope in $(r-4)_{21}\ $gives an element of $B_{r}$.Therefore, the
cardinality of $B_{r}\ $is at least the total number of $(r-1)$-crosspolytopes
in $(r-4)_{21}\ $which is $\left[  E_{r}:D_{r-1}\right]  \ $(see
section\ \ref{Section-Polytope}).

Now we want to show $\left\vert B_{r}\right\vert =\left[  E_{r}:D_{r-1}%
\right]  $ by calculating the number of pairs of lines with intersection $1$.

Since $W(S_{r})\ $acts transitively on $L_{r}$, we can focus on a line
$e_{r}\ $in $Pic\ S_{r}$ . The isotropy subgroup in $W(S_{r})\ $of $e_{r}\ $is
$E_{r-1}$-type as we saw in subsection \ref{SubSec-Constuction}. We choose
another line $h-e_{1}-e_{r}\ $with $\left(  h-e_{1}-e_{r}\right)  \cdot
e_{r}=1$.\ Since $W(S_{r})\ $preserves the intersection, the isotropy subgroup
$E_{r-1}\ $of $e_{r}\ $has an orbit of $h-e_{1}-e_{r}\ $such that each line in
the orbit has intersection $1\ $with $e_{r}$.\ The isotropy subgroup of
$h-e_{1}-e_{r}\ $in $E_{r-1}$, the isotropy subgroup of $e_{r}$, is generated
by the subgraph of Dynkin diagram which is perpendicular to $h-e_{1}-e_{r}$,
i.e.
\[%
\begin{tabular}
[c]{lllllll}
&  & $\underset{d_{0}}{{\Huge \cdot}}$ &  &  &  & \\
&  & $\ {\Huge \shortmid}$ &  &  &  & \\
$\underset{d_{2}}{{\Huge \cdot}}$ & ${\Huge -}$ & $\underset{d_{3}%
}{{\Huge \cdot}}$ & ${\Huge -}$ & ${\Huge ...}$ & ${\Huge -}$ & $\underset
{d_{r-2}}{{\Huge \cdot}}\ $.
\end{tabular}
\ \ \
\]
Therefore it is $D_{r-2}$-type,\ and more the total number of elements in the
orbit is given by $\left[  E_{r-1}:D_{r-2}\right]  $. In fact, all the lines
in $L_{r}\ $intersecting $e_{r}\ $by $1\ $is in this orbit by lemma
\ref{lemma-int2}. Thus, we obtain that the number of lines intersecting by
$1\ $with a fixed line is $\left[  E_{r-1}:D_{r-2}\right]  $.

Now, the number of pairs in $L_{r}\ $with intersection $1\ $is given by%
\begin{align*}
&  \text{ }\frac{1}{2}\left\vert L_{r}\right\vert \cdot\text{number of lines
intersecting by }1\ \text{with a fixed line}\\
&  =\frac{1}{2}\left[  E_{r}:E_{r-1}\right]  \cdot\left[  E_{r-1}%
:D_{r-2}\right]  =\frac{1}{2}\left[  E_{r}:D_{r-2}\right]  \text{.}%
\end{align*}
On the other hand, the total number of pairs of lines with intersection
$1\ $given by the $(r-1)$-crosspolytopes in $(r-4)_{21}\ $is obtained from
\begin{align*}
&  \text{ }\frac{1}{2}\left\vert \{(r-1)\text{-crosspolytope\ in }%
(r-4)_{21}\}\right\vert \\
&  \cdot\left\vert \{\ \text{pairs of lines with intersection }1\ \text{in an
}(r-1)\text{-crosspolytope}\}\right\vert \\
&  =\frac{1}{2}\left[  E_{r}:D_{r-1}\right]  \cdot\left[  D_{r-1}%
:D_{r-2}\right]  =\frac{1}{2}\left[  E_{r}:D_{r-2}\right]  \text{.}%
\end{align*}
Therefore, all the pairs of lines with intersection $1\ $in $Pic\ S_{r}\ $are
given from $(r-1)$-crosspolytopes in $(r-4)_{21}$. This gives
theorem.\ \ \ \ \ \ ${\LARGE \blacksquare}$

\bigskip

\section{The Geometry of Gosset Polytopes in $Pic\ S_{r}\otimes\mathbb{Q}$}

In this section, we explain the correspondences between subpolytopes of Gosset
polytopes $(r-4)_{21}$ and certain divisor classes in del Pezzo surfaces.

Here we consider the following divisor classes on a del Pezzo surface $S_{r}$
satisfying one of the following conditions.
\begin{align*}
\text{(1)\ }L_{r}^{a}  &  :=\{D\in Pic\ S_{r}\mid D\ =l_{1}+...+l_{a},\text{
}l_{i}\ \text{disjoint lines in }L_{r}\}\\
\text{(2)\ }\mathcal{E}_{r}  &  :=\{D\in Pic\ S_{r}\mid D^{2}=1,\text{
}K_{S_{r}}\cdot D=-3\}\\
\text{(3) }F_{r}  &  :=\{D\in Pic\ S_{r}\mid D^{2}=0,\text{ }K_{S_{r}}\cdot
D=-2\}
\end{align*}

We call the divisor classes in $L_{r}^{a}$ \textit{skew }$a$\textit{-lines}
for each $a$. Note that a skew $a$-line $D$ satisfies $D^{2}=-a\ $and
$K_{S_{r}}\cdot D=-a$. In fact, as we will see, skew $a$-lines are equivalent
to the divisors with $D^{2}=-a\ $and $K_{S_{r}}\cdot D=-a$, when $1\leq
a\leq3$.

The linear system of a divisor in the class $D\ $with $D^{2}=1$, $K_{S_{r}%
}\cdot D=-3$ induce a regular map to $\mathbb{P}^{2}\ $(\cite{Dolgachev}). We
call a divisor class in $Pic\ S_{r}$ with these conditions an
\textit{exceptional system}. As we will see, each choice of disjoint lines in
$S_{r}$ which producing blow-downs from $S_{r}$ to $\mathbb{P}^{2}$ gives one
of these divisor classes, and the converse is also true for $r<8$. Note that
for $S_{6}$, each effective divisor $\tilde{D}$ with $\tilde{D}^{2}%
=1\ $and~$K_{S_{r}}\cdot\tilde{D}=-3$ is corresponded to the twisted cubic
surface (see \cite{Buck-Kos}).

The divisor in class $D$ in $F_{r}\ $with $D^{2}=0,$ $K_{S_{r}}D=-2$ is
corresponded to the fiber class of a \textit{ruling} on $S_{r}$. Here a ruling
is a fibration of $S_{r}$ over $\mathbb{P}^{1}$ whose generic fiber is a
smooth rational curve. The rulings on the del Pezzo surfaces $S_{r}$ are
studied on behalf of researching on the line bundles on del Pezzo surfaces
according to representation theory (see \cite{Leung1}\cite{Leung-Zhang}). We
also call the divisor classes in $F_{r}$ \textit{rulings}.

\bigskip

\textbf{Dual lattices and theta functions}

The cardinalities of $L_{r}$, $L_{r}^{a}$, $\mathcal{E}_{r}$ and $F_{r}\ $are
finite because each set is an integral subset of a sphere in an affine
hyperplane section in $Pic\ S_{r}$. Furthermore, one can obtain them by
solving the corresponding Diophantine equations. But when $r\geq6$, the
cardinalities are very big so that we need another way to guarantee the
validity of our results of the simple calculation. The following argument
deals somewhat different objects with what we are aiming. Thus we sketch the
key ideas only and discuss the detail in another article. In the below
argument, we follow the notations in \cite{Conway-Sloane}.

First of all, we observe that the hyperplane section $H_{0}\ $in
$Pic\ S_{r}\ $is the root lattice $\Gamma_{r}\ $of the root system of $E_{r}$,
$r\geq6$. Here the set of minimal vectors in $\Gamma_{r}\ $is the set of
roots. We consider the translation of $H_{b}=\tilde{H}_{b}\cap Pic\ S_{r}$
along the vector $K_{S_{r}}\ $in $Pic\ S_{r}\otimes\mathbb{Q}$ .\ If $\frac
{b}{9-r}K_{S_{r}}\ $is an integral vector, then $H_{b}+\frac{b}{9-r}K_{S_{r}%
}\ $in $\tilde{H}_{0}\ $is $\Gamma_{r}$. When $r=7\ $and $b\ \in
\mathbb{Z-}2\mathbb{Z\ }$(resp. $r=6\ $and $b\ \in\mathbb{Z-}3\mathbb{Z}$),
there is an a glue vector $v\ $in $\tilde{H}_{0}\ $with norm $3/2\ $(resp.
$4/3$) such that $v+\Gamma_{r}\ $is $H_{b}+\frac{b}{9-r}K_{S_{r}}$. Therefore,
for each $b\ \in\mathbb{Z}$,\ $H_{b}\ $is transformed into the union of
$v+\Gamma_{r}\ $in $\tilde{H}_{0}\ $where $v\ $is a glue vector which can be
the origin. In fact, this union is the dual lattice of $\Gamma_{r}$. Now, in
order to figure out the total number of integral solutions $D\in Pic\ S_{r}%
\ $satisfying $D^{2}=a\ $and $D\cdot K_{S_{r}}=-b$, we send $D\ $to $\tilde
{H}_{0}\ $by $D+\frac{b}{9-r}K_{S_{r}}\ $and compute
\[
-\left(  D+\frac{b}{9-r}K_{S_{r}}\right)  ^{2}=-a+\frac{b^{2}}{9-r}\text{.}%
\]
And we look up the coefficient of degree $-a+\frac{b^{2}}{9-r}\ $in the theta
series of dual lattice of $\Gamma_{r}$. This gives us the number of solutions
for the equations.

For example, to identify the number of elements in $L_{r}$ via theta series,
we observe\
\begin{align*}
K_{S_{r}}l  &  =-1,\ l^{2}=-1\\
&  \Longleftrightarrow(l+\frac{K_{S_{r}}}{9-r})\cdot K_{S_{r}}=0,(l+\frac
{K_{S_{r}}}{9-r})^{2}=-1-\frac{1}{9-r},
\end{align*}
and the number of lines in $Pic\ S_{r}\ $is the same with the the coefficient
of degree $\left(  1+\frac{1}{9-r}\right)  $ in the theta series of the dual
lattices of the root lattice of $E_{r}$. If we consider $L_{7}$, the divisor
classes in $L_{7}\ $is transformed into the dual lattice of $\Gamma_{7}\ $by
$D+1/2K_{S_{7}}$ with norm $3/2$. And the coefficient of the degree $3/2\ $in
the theta series of the dual lattice of $E_{7}\ $is $56\ $%
(\cite{Conway-Sloane}).

Note that the root lattice $\Gamma_{8}\ $of $E_{8}\ $is self dual, and the
theta series of self dual lattice $\Gamma_{8}$ is given as%
\[
\Theta_{\Gamma_{8}}=%
%TCIMACRO{\dsum \limits_{m=0}^{\infty}}%
%BeginExpansion
{\displaystyle\sum\limits_{m=0}^{\infty}}
%EndExpansion
N_{m}q^{m},\ N_{m}=240\sigma_{3}(\frac{m}{2})
\]
where $\sigma_{r}(m)=%
%TCIMACRO{\dsum \nolimits_{d\mid m}}%
%BeginExpansion
{\displaystyle\sum\nolimits_{d\mid m}}
%EndExpansion
d^{r}$. As we will see, lines and rulings in $S_{8}$ correspond to lattice
points with norm $2$ and $4$ respectively. And the coefficients of $q^{2}$ and
$q^{4}$ are $240$ and $240(1+2^{3})=2160$.

\bigskip

\subsection{Special Divisor classes of del Pezzo surface $S_{r}\ $and
Subpolytopes of Gosset Polytope $(r-4)_{21}$}

\bigskip

In this subsection, we show $L_{r}^{a}$, $\mathcal{E}_{r}$ and $F_{r}\ $ are
bijective to subpolytopes in $(r-4)_{21}$. By lemmar \ref{LMAHyperp-Refl} ,
the conditions in $L_{r}^{a}$, $\mathcal{E}_{r}$ and $F_{r}$ are preserved by
the action of $W(S_{r})\ $the Weyl group of the root spaces in $Pic\ S_{r}$.
Therefore, the correspondences between these divisor classes on del Pezzo
surfaces and the subpolytopes in $(r-4)_{21}$ will be more than numerical
coincidences. Furthermore, since the geometry of subpolytopes in $E_{r}%
$-semiregular polytopes is basically configuration of vertices, it is natural
to study the divisor classes in $Pic\ S_{r}$ with respect to the
configurations of lines.\bigskip

There are two types of subpolytopes in $(r-4)_{21}$, which are $a\ $%
-simplexes\ $(0\leq a\leq r-1)\ $and ($r-1$)-crosspolytopes. In particular,
the facets of $(r-4)_{21}\ $consist of ($r-1$)-simplexes and ($r-1$%
)-crosspolytopes. In this subsection, these top degree subpolytopes
corresponds to the divisor classes representing rational maps from $S_{r}\ $to
$\mathbb{P}^{2}\ $and $\mathbb{P}^{1}\times$ $\mathbb{P}^{1}\ $respectively.

\bigskip

\textbf{Skew }$a$\textbf{-lines and }$(a-1)$\textbf{-simplexes}

In the subsection \ref{Sebsec-geomPoly}, $\ $we show that $A_{r}^{a}\ $the set
of centers of $a$-simplexes in $(r-4)_{21}\ $is bijective to the set of
$a$-simplexes in $(r-4)_{21}$. In fact the set $A_{r}^{a}\ $is the same with
the set of divisor classes so-called \textit{skew }$(a+1)$\textit{-lines},
\[
L_{r}^{a+1}:=\{D\in Pic(S_{r})\mid D=l_{1}+...+l_{a+1},\text{ }l_{i}%
\ \text{disjoint lines in }L_{r}\}.
\]
Therefore, we have the following theorem.

\bigskip

\begin{theorem}
\label{alinevssimplex}The set of skew $a$-lines in $Pic\ S_{r}$, $1\leq a\leq
r$ is bijective to the set of $\left(  a-1\right)  $-simplexes, $\alpha
_{(a-1)}$, in the Gosset polytope $(r-4)_{21}$.
\end{theorem}

\bigskip

\textbf{Remark} : \textbf{1}.\ In particular, a skew $1$-line is just a $line$
in $L_{r}$ which correspond to the vertices of $E_{r}$-semiregular polytope,
and skew $2$-lines represent edges in $(r-4)_{21}$.

\textbf{2}. Since the vertices of an $\left(  a-1\right)  $-simplex correspond
to the lines consisting of a skew $a$-line in $S_{r}$, each $\left(
a-1\right)  $-simplex represents a rational map from $S_{r}$ to $S_{r-a}$
obtained by blowing down the disjoint $a$-lines in $S_{r}$. In particular,
each $\left(  r-1\right)  $-simplex in $(r-4)_{21}$ gives a rational map from
$S_{r}$ to $\mathbb{P}^{2}$.

\textbf{3}. As the Weyl group $W(S_{r})\ $act transitively on the set of
$\left(  a-1\right)  $-simplexes in the Gosset polytope $(r-4)_{21}$,
$W(S_{r})\ $acts transitively on the set of skew $a$-lines.

\bigskip

It is easy to see that skew $a$-lines in $Pic\ S_{r}\ $satisfies $D^{2}%
=-a\ $and $K_{S_{r}}D=-a$, which is called $a$\textit{-divisors}. The number
of $a$-divisors in $S_{r},$ $2\leq a\leq r,$ can be obtained from the theta
series of dual lattice of the root lattice of $E_{r}$. First we observe%
\begin{align*}
K_{S_{r}}D  &  =-a,\ D^{2}=-a\\
&  \Longleftrightarrow(D+a\frac{K_{S_{r}}}{9-r})\cdot K_{S_{r}}=0,(D+a\frac
{K_{S_{r}}}{9-r})^{2}=-a\left(  1+\frac{a}{9-r}\right)  ,
\end{align*}
and according to \cite{Conway-Sloane}, we have the table of the number of
$2$-divisors and $3$-divisors.%
\[
\underset{\text{Number of }2\text{-and }3\text{-divisors in }S_{r}}{%
\begin{tabular}
[c]{|l||l|l|l|l|l|l|}\hline
$a$ & $S_{3}$ & $S_{4}$ & $S_{5}$ & $S_{6}$ & $S_{7}$ & $S_{8}$\\\hline
$2$ & $6$ & $30$ & $80$ & $216$ & $756$ & $6720$\\\hline
$3$ & $2$ & $30$ & $160$ & $720$ & $4032$ & $60480$\\\hline
\end{tabular}
\ }%
\]
In fact, the table is parallel to the number of $1$-simplexes and
$2$-simplexes in Gosset polytopes $(r-4)_{21}\ $(see section
\ref{Section-Polytope}). Therefore, by above lemma and theorem, the divisor
classes $D$ in $Pic\ S_{r}$ with $D^{2}=-a\ $and $K_{S_{r}}D=-a$, $1\leq
a\leq3$ are skew $a$-lines. This gives the following theorem.

\bigskip

\begin{theorem}
For $1\leq a\leq3$, each divisor class $D$ in $Pic\ S_{r}$ with $D^{2}%
=-a\ $and $K_{S_{r}}D=-a$ can be written as the sum of skew lines
$l_{1},..,l_{a}$ where the choice is unique up to the permutation.
\end{theorem}

\bigskip

\textbf{Exceptional systems and }$(r-1)$\textbf{-simplexes in }$(r-4)_{21}$

Recall $\left(  r-1\right)  $-simplexes in $(r-4)_{21}\ $are one of the two
types facets in $(r-4)_{21}$, and each $\left(  r-1\right)  $-simplex in
$(r-4)_{21}$ corresponds to disjoint $r$-lines in $Pic\ S_{r}\ $gives a
rational map from $S_{r}$ to $\mathbb{P}^{2}$. But the above correspondence
between $\left(  r-1\right)  $-simplex in $(r-4)_{21}\ $and skew $r$-lines in
$Pic\ S_{r}\ $is somewhat coarse and we want to have another approach. Here we
consider classes which we called exceptional systems.

An \textit{exceptional system} is a divisor class $D$ on a del Pezzo surface
$S_{r}$ with $D^{2}=1$ and $D\cdot K_{S_{r}}=-3$, and the linear system of $D$
gives a regular map from $S_{r}$ to $\mathbb{P}^{2}$ (\cite{Dolgachev}). We
denote the set of exceptional systems in $Pic\ S_{r}$ as
\[
\mathcal{E}_{r}:=\{D\in Pic(S_{r})\mid D^{2}=1,K_{S_{r}}\cdot D=-3\}.
\]
When $r=6$, each linear system with above conditions contains a twisted cubic
curve and it is also known there are $72$ of such classes. In fact, according
to the correspondence between skew $a$-lines in $S_{6}$ and $(a-1)$--simplexes
in $2_{21}$, we can see the number of exceptional systems in $S_{6}$ equals to
the number of $5$-simplexes. As we saw in the remark of theorem
\ref{alinevssimplex}, each $(r-1)$-simplex in $(r-4)_{21}\ $corresponds to the
rational map from $S_{r}\ $to $\mathbb{P}^{2}$, therefore it is natural to
compare the set of $(r-1)$-simplex in $(r-4)_{21}\ $and the set of exceptional
systems. Here, we explain each $(r-1)$-simplex in $(r-4)_{21}$ is related to
an exceptional system, and more two sets of these are bijective for $3\leq
r\leq7$.

\bigskip

First, we observe
\begin{align*}
D\cdot K_{S_{r}}  &  =-3\text{, }D^{2}=1\text{ }\\
&  \Longleftrightarrow\left(  D+\frac{3K_{S_{r}}}{9-r}\right)  \cdot K_{S_{r}%
}=0,\left(  D+\frac{3K_{S_{r}}}{9-r}\right)  ^{2}=1-\frac{9}{9-r}\text{,}%
\end{align*}
and from the theta series of dual lattice of root lattice of $E_{r}$, the
number of exceptional systems can be listed as below. We observe that the
numbers of exceptional systems in del Pezzo surfaces and the numbers of top
degree subsimplexes in $(r-4)_{21}$ are parallel except $r=8$.%

\[%
\begin{tabular}
[c]{|l|l|l|l|l|l|l|}\hline
del Pezzo Surfaces $S_{r}$ & $\ \ S_{3}$ & $S_{4}$ & $S_{5}$ & $S_{6}$ &
$S_{7}$ & $S_{8}$\\\hline
number of exceptional systems & $\ \ 2$ & $5$ & $16$ & $72$ & $576$ &
$17520$\\\hline\hline
$(r-4)_{21}$ & $-1_{21}$ & $0_{21}$ & $1_{21}$ & $2_{21}$ & $3_{21}$ &
$4_{21}$\\\hline
number of $(r-1)$-simplexes & $\ \ 2$ & $5$ & $16$ & $72$ & $576$ &
$17280$\\\hline
\end{tabular}
\]

To explain the correspondence, we define a transformation $\Phi\ $from
$\mathcal{E}_{r}\ $to $L_{r}^{r}\ $by
\[
\Phi(D_{t})=:K_{S_{r}}+3D_{t}\ \text{for }D_{t}\in\mathcal{E}_{r}\text{.}%
\]
It is well-defined because
\[
\Phi(D_{t})\cdot K_{S_{r}}=\left(  K_{S_{r}}+3D_{t}\right)  \cdot K_{S_{r}%
}=-r\ \text{,\ }\Phi(D_{t})^{2}=-r\text{.}%
\]
In the following theorem, the transformation $\Phi\ $leads a correspondence
between $\mathcal{E}_{r}\ $and $L_{r}^{r}$.

\begin{theorem}
When $3\leq r\leq$ $8$, each $(r-1)$-simplex in $(r-4)_{21}$ corresponds to an
exceptional system in the del Pezzo surfaces $S_{r}$. Moreover, for $3\leq
r\leq$ $7$, the Weyl group $W(S_{r})$ act transitively on $\mathcal{E}_{r}$
the set of exceptional systems in the del Pezzo surface $S_{r}$, and
$\mathcal{E}_{r}$ is bijective to $L_{r}^{r}$ the set of skew $r$-lines in
$Pic$ $S_{r}$.
\end{theorem}

\textbf{Proof}: First, we observe the above transformation $\Phi\ $is
injective and equivariant for the action of the Weyl group $W(S_{r})\ $by
lemma\ref{LMAHyperp-Refl}. We consider a class $h$ in $Pic(S_{r})\ $which is
an exceptional system. Here $\Phi\ $sent $h$ to the sum of disjoint
exceptional classes $e_{1}+...+e_{r}$ which is a skew $r$-line representing
one of $(r-1)$-simplexes by theorem \ref{alinevssimplex}. And the orbit
containing $h\ $in $\mathcal{E}_{r}\ $is corresponded to the orbit containing
$e_{1}+...+e_{r}\ $in $L_{r}^{r}$. Since the action of $W(S_{r})\ $on
$L_{r}^{r}\ $is transitive, each skew $r$-line representing an $(r-1)$-simplex
corresponds to an exceptional system. According to the above table,
$\mathcal{E}_{r}\ $is also transitively acted $W(S_{r})\ $for $3\leq r\leq$
$7$. ${\LARGE \blacksquare}$

\bigskip

\textbf{Remark:} $\ $When $r=8$, the set of exceptional systems has two
orbits. One orbit with $17280$ elements corresponds to the set of skew
$8$-lines in $S_{8}$, and the other orbit with $240$ elements corresponds to
the set of $E_{8}$-roots since for each $E_{8}$-root $d$, $-3K_{S_{8}}+$ $d$
is an exceptional system.

\bigskip

\textbf{Rulings and Crosspolytopes}

\bigskip

Now, the crosspolytopes which is the other type of facets in $(r-4)_{21}\ $is
only remained from the subpolytopes in $(r-4)_{21}$.\ In subsection
\ref{Sebsec-geomPoly}, the set of $(r-1)$-crosspolytopes in $(r-4)_{21}\ $is
bijective to the set of their centers defined by $B_{r}$. Observe that each
element\ in $B_{r}$, $l_{1}+l_{2}\ $with $l_{1}\cdot l_{2}=1$, satisfies
\[
\left(  l_{1}+l_{2}\right)  ^{2}=0\ ,K_{S_{r}}\cdot\left(  l_{1}+l_{2}\right)
=-2\text{,}%
\]
and we consider the divisor classes with these conditions. Furthermore, the
divisor classes correspond to the rational maps from $S_{r}\ $to
$\mathbb{P}^{1}\times$ $\mathbb{P}^{1}$.

The divisor class $f$ on del Pezzo surface $Pic\ S_{r}$ with $f^{2}=0,$
$K_{S_{r}}\cdot f=-2$ is called a \textit{ruling} since the divisor in it
corresponds to a fibration of $S_{r}$ over $\mathbb{P}^{1}$ whose generic
fiber is a smooth rational curve. The set of rulings in $Pic\ S_{r}$ is
denoted as
\[
F_{r}:=\{f\in Pic(S_{r})\mid f^{2}=0,K_{S_{r}}\cdot f=-2\}.
\]

The number of rulings in $Pic\ S_{r}$ can be obtained from the theta series of
the dual lattice of root lattice of $E_{r}$ as followings
\begin{align*}
K_{S_{r}}f  &  =-2,\ f^{2}=0\\
&  \Longleftrightarrow\left(  f+\frac{2K_{S_{r}}}{9-r}\right)  \cdot K_{S_{r}%
}=0,\ \left(  f+\frac{2K_{S_{r}}}{9-r}\right)  ^{2}=-\frac{4}{9-r}.
\end{align*}
Furthermore, we get the following parallel list of the numbers of rulings in
del Pezzo surface $S_{r}$ and the numbers of crosspolytopes in $(r-4)_{21}$.%
\[%
\begin{tabular}
[c]{|l|l|l|l|l|l|l|}\hline
del Pezzo Surfaces $S_{r}$ & \ $S_{3}$ & $S_{4}$ & $S_{5}$ & $S_{6}$ & $S_{7}$
& $S_{8}$\\\hline
number of rulings & $3$ & $5$ & $10$ & $27$ & $126$ & $2160$\\\hline\hline
$(r-4)_{21}$ & $-1_{21}$ & $0_{21}$ & $1_{21}$ & $2_{21}$ & $3_{21}$ &
$4_{21}$\\\hline
number of $\left(  r-1\right)  $-crosspolytopes & $3$ & $5$ & $10$ & $27$ &
$126$ & $2160$\\\hline
\end{tabular}
\]

To get a rough idea about these parallel lists of numbers, we consider the
followings. By lemma \ref{LMAHyperp-Refl}, the set of rulings on $S_{r}$ is
acted by the Weyl group $W(S_{r})$. We consider a ruling $h-e_{1}$. And the
isotropy group of $h-e_{1}$ is generated by the simple roots perpendicular to
$h-e_{1}$ with the relationships represented as
\[%
\begin{tabular}
[c]{lllllllll}
&  &  &  & $\underset{d_{0}}{{\Huge \cdot}}$ &  &  &  & \\
&  &  &  & $\ {\Huge \shortmid}$ &  &  &  & \\
$\ast$ & ${\large \cdots}$ & $\underset{d_{2}}{{\Huge \cdot}}$ & ${\Huge -}$ &
$\underset{d_{3}}{{\Huge \cdot}}$ & ${\Huge -}$ & ${\Huge ...}$ & ${\Huge -}$
& $\underset{d_{r-1}}{{\Huge \cdot}}$%
\end{tabular}
\ \
\]
Thus it is $D_{r-1}$-type, and the number of element in the orbit containing
$h-e_{1}$ is $\left[  E_{r}:D_{r-1}\right]  $. This procedure is exactly same
with counting the number of $\left(  r-1\right)  $-crosspolytopes in an
$(r-4)_{21}$. As a matter of fact the action of the Weyl group $W(S_{r})\ $on
rulings is transitive from the following theorem. A part of following theorem
appears in \cite{Leung1} and \cite{Batyrev-Popov}.

\begin{theorem}
\label{Thm-ruling}For each ruling $f$ in a del Pezzo surface $S_{r}$, there is
a pair of lines $l_{1}$ and $l_{2}$ with $l_{1}\cdot l_{2}=1\ $such that $f$
is equivalent to the sum of $l_{1}+l_{2}$. Furthermore, the set of rulings in
$S_{r}$ is bijective to the set of $\left(  r-1\right)  $-crosspolytopes in
$(r-4)_{21}\ $and acted transitively\ by the Weyl group $W(S_{r})$.
\end{theorem}

\textbf{Proof}:\ From the subsection \ref{Sebsec-geomPoly} , $B_{r}$ (the set
of centers of $\left(  r-1\right)  $-crosspolytopes in $(r-4)_{21}$) is
already a subset of rulings $F_{r}$. Since $B_{r}\ $is bijective to the set of
$\left(  r-1\right)  $-crosspolytopes in $(r-4)_{21}$,\ $\left\vert
B_{r}\right\vert =\left\vert F_{r}\right\vert $. Therefore $B_{r}=F_{r}$ and
each ruling can be written as the two lines with intersection $1$.\ This also
show that $F_{r}\ $is bijective to the set of $\left(  r-1\right)
$-crosspolytopes in $(r-4)_{21}$. Since this correspondence is naturally
equivariant for the Weyl group $W(S_{r})\ $action and the set of $\left(
r-1\right)  $-crosspolytopes in $(r-4)_{21}\ $acted transitively by $W(S_{r}%
)$, $F_{r}\ $is also acted transitively by $W(S_{r})$.
\ \ ${\LARGE \blacksquare}$

\bigskip

According to the theorem \ref{Correspondence line and vetices}, the vertices
of an $\left(  r-1\right)  $-crosspolytope must correspond lines in the linear
system of a ruling. Here we consider the pairs of antipodal vertices in the
$\left(  r-1\right)  $-crosspolytope and their correspondences in a ruling.
Each antipodal vertices in the pair are the two bipolar points in the
crosspolytope which is a bipyramid. And the number of the pairs of antipodal
vertices in an $\left(  r-1\right)  $-crosspolytope is $\left(  r-1\right)  $.
Therefrom, we get the following corollary.

\begin{corollary}
\label{Coro-ruling-crosspl}For each ruling $f$ in a del Pezzo surface $S_{r}$,
there are $(r-1)$-pairs of lines with $1$-intersection whose sum is $f$.
Furthermore, each of these pairs corresponds to antipodal vertices of the
$(r-1)$-crosspolytope corresponding to $f$ in $(r-4)_{21}$.
\end{corollary}

\bigskip

From the above corollary, we obtain the following useful lemma.

\bigskip

\begin{lemma}
\label{lemma-ruling}For a ruling $f$ and a line $l$ in $Pic$ $S_{r}$, the
following are equivalent.

(1) $f\cdot l=0$

(2) $f-l$ is a line

(3) The vertex$\ $represented by $l$ in $(r-4)_{21}$\ is one of the vertex of
the $(r-1)$-crosspolytope corresponding to $f$.
\end{lemma}

\bigskip

\textbf{Remark}: As the top degree simplexes in $(r-4)_{21}$ are related to
the blowing down from $S_{r}$ to $\mathbb{P}^{2}$, the crosspolytope which is
the other type of facets in $(r-4)_{21}$ gives the blowing down from $S_{r}$
to $\mathbb{P}^{1}\times\mathbb{P}^{1}$. Here we give an example of this blow
down and the further study will be discussed in \cite{Clinger-Lee}.

For cubic surfaces $S_{6}$, it is well known that there are two disjoint lines
$l_{a}$ and $l_{b}$, and five parallel lines $l_{i},$ $1\leq$ $i\leq5$ meet
both $l_{a}$ and $l_{b}$. Furthermore, for each pair $l_{a}$ (resp. $l_{b}$)
and $l_{i}$, there is a line $l_{ai}$ (resp. $l_{bi}$) intersecting $l_{a}$
(resp. $l_{b}$) and $l_{i}$. Blowing down of five disjoint lines $l_{i}$ gives
a rational map from $S_{6}$ to $\mathbb{P}^{1}\times\mathbb{P}^{1}$. From the
facts in this article, for each line $l_{i}$, $l_{ai}+l_{i}$ produces the same
ruling, namely, $l_{i}$ and $l_{ai}$ for $1\leq$ $i\leq5$ correspond to ten
vertices of a $5$-crosspolytope. And the same fact is true for lines $l_{i}$
and $l_{bi}$, $1\leq$ $i\leq5$. Therefore proper choice of $5$-crosspolytopes
in a Gosset polytope $2_{21}$ gives a blowing down from $S_{6}$ to
$\mathbb{P}^{1}\times\mathbb{P}^{1}$.

For example, we consider two disjoint lines $e_{6}$ and $2h-\sum_{k=1}%
^{6}e_{k}+e_{6}$ in $S_{6}$. For the line $e_{6}$, there are $5$-pairs of
lines $e_{i}$ and $h-e_{i}-e_{6},\ 1\leq i\leq5\ $where $e_{6}$ , $e_{i}$ and
$h-e_{i}-e_{6}$ have $1$-intersections to each other. Similarly, the line
$2h-\sum_{k=1}^{6}e_{k}+e_{i}$ has $5$-pairs of lines $2h-\sum_{k=1}^{6}%
e_{k}+e_{i}$ and $h-e_{i}-e_{6},\ 1\leq i\leq5$. And we find $5$-disjoint
lines $h-e_{i}-e_{6},\ 1\leq i\leq5$ which give a blowing down to
$\mathbb{P}^{1}\times\mathbb{P}^{1}$. Here we observe the five pair of lines
for $e_{6}$ are from a ruling $\left(  h-e_{i}-e_{6}\right)  +e_{i}=h-e_{6}$
and corresponded to the antipodal vertices in the crosspolytope in $2_{21}$
corresponding to $h-e_{6}$. Similarly the five pair of lines for
$h-e_{i}-e_{6}$ are from a ruling $\left(  h-e_{i}-e_{6}\right)  +\left(
2h-\sum_{k=1}^{6}e_{k}+e_{i}\right)  =3h-\sum_{k=1}^{6}e_{k}-e_{6}$.

This observation induces interesting relationships between geometry of cubic
surfaces and the combinatorial data on a Gosset polytope $2_{21}$.
Furthermore, as each del Pezzo surface blows down to $\mathbb{P}^{1}%
\times\mathbb{P}^{1}$, we can search the similar works on all the del Pezzo
surfaces and the corresponding Gosset polytopes. (see \cite{Clinger-Lee})

\bigskip

\bigskip

\section{Applications}

\subsection{Monoidal transforms for lines}

In this subsection, we consider lines along the blow-up procedure producing
del Pezzo surfaces to study the geometry of Gosset polytopes $(r-4)_{21}$
according to the above correspondences.

This blow-up procedure can be applied to rulings and other special divisor
classes in this article so as to obtain recursive relationships corresponded
to subpolytopes in Gosset polytopes $(r-4)_{21}$. Furthermore, we can also use
this procedure by virtue of interesting combinatorial relationships on
subpolytopes in Gosset polytopes $(r-4)_{21}$, which are also related to the
Cox ring of del Pezzo surfaces (\cite{Batyrev-Popov}). This will be explained
in \cite{Lee}.

For a fixed vertex $P\ $in a Gosset polytope $(r-4)_{21}$, the set of vertices
with the shortest distance from $P\ $is characterized as the vertex figure of
$P\ $by the action of an isotropy group of $E_{r-1}$-type. But the description
of the set of vertices with further distance from $P\ $requires rather
indirect procedure. Here, we consider divisor classes producing lines by
blow-up and apply these to study the set of vertices with further distance
from $P\ $according to the correspondences between vertices and lines.

For a fixed line $l$ in a del Pezzo surface $S_{r}$, $3\leq r\leq8$, we
consider a set
\[
N_{k}(l,S_{r})=\left\{  \ l^{\prime}\ \in L_{r}\mid l^{\prime}\cdot
l=k\right\}
\]
which also$\ $presents the set of lines in $S_{r}\ $with the same distance
from $l$. In section \ref{Sec-GosinPic}, we can describe local geometry of the
polytope via case study on these sets of lines. Here, by using blow-up
procedure on lines,we complete the description in the uniform manner.

Since a del Pezzo surface $S_{r}$ is obtained by blowing up one point on
$S_{r-1}\ $to a line $l\ $in $S_{r}$, we can describe divisor classes
in$\ S_{r-1}$ producing lines in $S_{r}$ after blowing up. In fact, the choice
of line $l$ in above can be replaced by exceptional class $e_{r}$ in $S_{r}$
which is the exceptional divisor in $S_{r}$ given by a blow-up of a point from
$S_{r-1}$. The proper transform of a divisor $D$ in $Pic(S_{r-1})$ producing a
line in $S_{r}\ $satisfies
\[
(D-me_{r})^{2}=-1\text{, }(D-me_{r})\cdot(K_{S_{r-1}}+e_{r})=-1
\]
for a nonnegative integer $m$. Therefore, we consider a divisor $D$ in
$Pic(S_{r-1})$ with
\[
D^{2}=m^{2}-1\text{, }D\cdot K_{S_{r-1}}=-m-1.
\]
By the Hodge index theorem, we have%
\[
(m^{2}-1)K_{S_{r-1}}^{2}=D^{2}K_{S_{r-1}}^{2}\leq(D\cdot K_{S_{r-1}}%
)^{2}=(-m-1)^{2},
\]
which implies%
\[
-1\leq m\leq1+\frac{2}{9-r}.
\]
Thus the list of possible $m$ is
\[
m=\left\{
\begin{tabular}
[c]{l}%
$0,1\ $\ $\ \ \ \ \ \ \ \ $\\
$0,1,2\ \ \ \ $\\
$0,1,2,3\ $%
\end{tabular}%
\begin{tabular}
[c]{l}%
$4\leq r\leq6$\\
$r=7$\\
$r=8$%
\end{tabular}
\ \ \right.  .
\]

\begin{definition}
\bigskip A line $l$ in the Picard group of del Pezzo surface $S_{r}$ obtained
by blow-up from a divisor class $D$ in $S_{r-1}$ is called an $m$%
\textit{-degree line} if $l$ $=D-me_{r}\ $where $e_{r}\ $is the exceptional
class produced by the blow-up.
\end{definition}

\bigskip

(1) $0$\textbf{-degree line in }$S_{r}$, $4\leq r\leq8\ $: Each $0$-degree
line in $S_{r}$ corresponds to a line in $S_{r-1}$ and the number of
$0$-degree lines equals to the number of lines in $S_{r-1}$.

(2) $1$\textbf{-degree line in }$S_{r},$ $4\leq r\leq8\ $: A divisor $D$ in
$S_{r-1}$ with $D^{2}=0\ $and $D\cdot K_{S_{r-1}}=-2$ corresponds to an
$1$-degree line in $S_{r}$. Therefore, $F_{r-1}$ the set of rulings in
$S_{r-1}$ is bijective to the set of $1$-degree lines in $S_{r}$.

(3) $2$\textbf{-degree line in }$S_{r},$ $r=7,8\ $:\ When $r=7$ (resp. $r=8$),
there are $2$-degree lines given by divisors in $Pic\ S_{6}$\ (resp.
$Pic\ S_{7}$) with $D^{2}=3\ $and $DK_{S_{6}}=-3$. For $r=7$, by the Hodge
index theorem, $-K_{S_{6}}$ is the only one divisor in $S_{6}$ with these
equations. For $r=8$, we can transform the equations to
\[
\left(  D+K_{S_{7}}\right)  ^{2}=-1,\ \left(  D+K_{S_{7}}\right)  \cdot
K_{S_{7}}=-1,
\]
which represent lines in $S_{7}$. Therefore, the number of $2$-degree line in
$S_{8}$ is the same with the number of lines in $S_{7}$.

(4) $3$\textbf{-degree line in }$S_{8}\ $:\ For $r=8$, there are $3$-degree
lines obtained from divisors in $Pic(S_{7})$ with $D^{2}=8\ $and $DK_{S_{7}%
}=-4$. These equations are equivalent to%
\[
\left(  D+2K_{S_{7}}\right)  ^{2}=0,\ \left(  D+2K_{S_{7}}\right)  \cdot
K_{S_{7}}=0,
\]
Therefore, $-2K_{S_{7}}$ is the only divisor in $S_{7}$ producing a $3$-degree
line in $S_{8}$.

\bigskip

Now, we obtain the following theorem.

\begin{theorem}
Let $l$ be a fixed line in a del Pezzo surface $S_{r}$ , $4\leq r\leq8$, and
$V_{l}$ be the vertex corresponding to the line $l$ in the polytope
$(r-4)_{21}$.

(1) $N_{0}(l,S_{r}),4\leq r\leq8$, is bijective to the set of lines $L_{r-1}$
in $S_{r-1},$ and equivalently, it is also bijective to the set of vertices in
the polytope $(r-5)_{21}$.

(2) $N_{1}(l,S_{r}),4\leq r\leq8$, is bijective to the set of ruling
containing $l\ $in $S_{r}$ which is also bijective to the set of rulings
$F_{r-1}$ in $S_{r-1}$. Equivalently, it is also bijective to the set of
$\left(  r-1\right)  $-crosspolytopes in the polytope $(r-4)_{21}$ containing
$V_{l}$ and the set of $\left(  r-2\right)  $-crosspolytopes in the polytope
$(r-5)_{21}$.

(3) $N_{2}(l,S_{8})$ is bijective to $N_{0}(-2K_{S_{8}}-l,S_{8})$ the set of
skew lines in $S_{7}$ for a line $-2K_{S_{8}}-l,$ and equivalently it is also
bijective to the set of lines in $S_{7}$.

(4) $N_{2}(l,S_{7})=\{-K_{S_{7}}-l\}\ $and $N_{3}(l,S_{8})=\{-2K_{S_{8}}-l\}$
\end{theorem}

\bigskip\textbf{Proof}: Since $N_{k}(l,S_{r})$ is equivalent to the $k$-degree
lines in $S_{r}$, the above description of blow-up for lines and theorem
\ref{Thm-ruling} give the theorem. ${\LARGE \blacksquare}$

\bigskip

\textbf{Remark}: 1.For $r=3$, $\left\vert N_{0}(l,S_{3})\right\vert =3$ and
$\left\vert N_{1}(l,S_{3})\right\vert =2$.

2. For $S_{7}$, each element $l^{\prime}$ in $N_{1}(l,S_{7})$ satisfies
$\left(  -K_{S_{7}}-l\right)  \cdot l^{\prime}=0$. Therefore, $N_{1}(l,S_{7})$
is bijective to $N_{0}(-K_{S_{7}}-l,S_{7})$, and the theorem explains
$\left\vert F_{6}\right\vert =\left\vert N_{1}(l,S_{7})\right\vert =\left\vert
N_{0}(-K_{S_{7}}-l,S_{7})\right\vert =$ $\left\vert L_{6}\right\vert $, namely
the number of lines and rulings in $S_{6}$ are the same.

3. For $S_{8}$, each line $l^{\prime}$ in $S_{8}$ satisfies $\left(
-2K_{S_{7}}-l\right)  \cdot l^{\prime}=2-l\cdot l^{\prime}$. Thus we have
$\left\vert N_{0}(l,S_{8})\right\vert =\left\vert N_{2}(-2K_{S_{8}}%
-l,S_{8})\right\vert $ and $\left\vert N_{1}(l,S_{8})\right\vert =\left\vert
N_{1}(-2K_{S_{8}}-l,S_{8})\right\vert $.

\bigskip

\subsection{Gieser Transform and Bertini Transform for lines}

Recall that for two lines $l_{1}$ and $l_{2}$ in $S_{7}$, $l_{1}\cdot l_{2}=2$
is equivalent to $l_{1}+l_{2}=-K_{S_{7}}$. In fact, these two lines in $S_{7}$
with $2$-intersection represent a bitangent of degree $2$ covering from
$S_{7}$ to $\mathbb{P}^{2}$ given by $\left\vert -K_{S_{7}}\right\vert $.
Similarly, any two lines $l_{1}^{\prime}$ and $l_{2}^{\prime}$ in $S_{8}$ with
$l_{1}^{\prime}\cdot l_{2}^{\prime}=3$ equivalently hold $l_{1}+l_{2}%
=-2K_{S_{8}}$ and these two lines in $S_{8}$ with $3$-intersection are related
to a tritangent plane of degree $2$ covering from $S_{8}$ to $\mathbb{P}^{2}$
given by $\left\vert -2K_{S_{8}}\right\vert $. Here the deck transformation
$\sigma$ of the cover is a biregular automorphism of $S_{7}$ (resp. $S_{8}$).
For a blow down $\pi:S_{7}\rightarrow\mathbb{P}^{2}$, $\pi\sigma$ is another
blow down, and the corresponding birational transformation $\left(  \pi
\sigma\right)  \pi^{-1}:\mathbb{P}^{2}\rightarrow\mathbb{P}^{2}$ is called
Gieser transform (resp. Bertini transform) (see chapter 8 \cite{Dolgachev}).
Therefore, Gieser transform corresponds to a transformation $G\ $on lines in
$S_{7}$ defined as%
\[
G(l):=-\left(  K_{S_{7}}+l\right)  ,
\]
and we also call $G\ $the \textit{Gieser transform on lines} or simply Gieser
transform. Similarly a transformation $B$ on lines in $S_{8}$ defined as
\[
B(l):=-\left(  2K_{S_{8}}+l\right)
\]
is referred as the \textit{Bertini transform on lines} or simply Bertini transform.

Since both Gieser transform and Bertini transform are defined on the set of
lines, we can extend the definition to any divisor written as linear sum of
lines. Namely, for a divisor class $D$ given as $a_{1}l_{1}+...+a_{m}l_{m}$ in
$S_{7}$%
\[
G(D):=a_{1}G\left(  l_{1}\right)  +...+a_{m}G\left(  l_{m}\right)  ,
\]
and similarly in $S_{8},$%
\[
B(D):=a_{1}B\left(  l_{1}\right)  +...+a_{m}B\left(  l_{m}\right)  .
\]
And for lines $l_{1},$ $l_{2}$ in $Pic\ S_{7}\ $and lines $l_{1}^{\prime}$,
$l_{2}^{\prime}\ $in $Pic\ S_{8}$, we have
\begin{align*}
G(l_{1})\cdot G(l_{2})  &  =\left(  K_{S_{7}}+l_{1}\right)  \cdot\left(
K_{S_{7}}+l_{2}\right) \\
&  =K_{S_{7}}^{2}+(l_{1}+l_{2})\cdot K_{S_{7}}+l_{1}\cdot l_{2}=l_{1}\cdot
l_{2}\text{,}\\
B(l_{1}^{\prime})\cdot B(l_{2}^{\prime})  &  =\left(  2K_{S_{8}}+l_{1}%
^{\prime}\right)  \cdot\left(  2K_{S_{8}}+l_{2}^{\prime}\right) \\
&  =4K_{S_{8}}^{2}+2(l_{1}^{\prime}+l_{2}^{\prime})\cdot K_{S_{8}}%
+l_{1}^{\prime}\cdot l_{2}^{\prime}=l_{1}^{\prime}\cdot l_{2}^{\prime}\text{.}%
\end{align*}
Therefore, $G$ and $B\ $preserves the intersections between lines, $G$ and $B$
are symmetries on $3_{21}$ and $4_{21}$, respectively. Moreover, since all the
regular subpolytopes we discussed in $3_{21}$ and $4_{21}$ are written as
linear sums of lines, $G$ and $B$ act on the set of these subpolytopes. In
summary, we have the following theorem.

\begin{theorem}
The Gieser transform $G\ $on the set lines in $Pic\ S_{7}$ and Bertini
transform $B$ on the set of lines in $Pic$ $S_{8}$ can be extended to a
symmetry of $3_{21}$ and $4_{21}$ respectively. Furthermore, $G$ and $B$ acts
on the each set of subsimplexes and the set of $6$-crosspolytopes (
$7$-crosspolytopes, respectively).
\end{theorem}

\bigskip

Naturally, the further studies on Gieser transform $G$ and Bertini
transform\ $B\ $are performed along the degree $2$ covering from $S_{7}$ to
$\mathbb{P}^{2}$ and $2$ covering from $S_{8}$ to $\mathbb{P}^{2}$. This will
be continued in \cite{Clinger-Lee}\cite{Lee}.

\bigskip

\textbf{Remark} : The above pairs of lines are special cases of Steiner blocks
which are related to the inscribed simplexes in $(r-4)_{21}$. This relation
will be discussed in \cite{Lee}.

\bigskip

\textit{Acknowledgments: Author is always grateful to his mentor Naichung
Conan Leung for everything he give to the author. Author expresses his
gratitude to Adrian Clingher for his kind helps, and thanks to Prabhakar Rao
and Ravindra Girivaru for useful discussions.}

\bigskip

\bigskip

{\scriptsize Addresses: }

{\scriptsize Jae-Hyouk Lee (lee@math.umsl.edu)}

{\scriptsize Department of Mathematics and Computer Science, University of
Missouri-St. Louis, U.S.A.}

\end{document}